\documentclass[a4paper,11pt]{article}
\title{\bf{Stability conditions and 
Calabi-Yau fibrations}
}
\date{}
\author{Yukinobu Toda}
\usepackage{amsmath,amssymb}
\usepackage{macros}

\usepackage{array}
\usepackage{amscd}
\usepackage[all]{xy}

\DeclareFontFamily{U}{rsfs}{%
\skewchar\font127}
\DeclareFontShape{U}{rsfs}{m}{n}{%
<-6>rsfs5<6-8.5>rsfs7<8.5->rsfs10}{}
\DeclareSymbolFont{rsfs}{U}{rsfs}{m}{n}
\DeclareSymbolFontAlphabet
{\mathrsfs}{rsfs}
\DeclareRobustCommand*\rsfs{%
\@fontswitch\relax\mathrsfs}
\newtheorem{thm}{Theorem}[section]
\newtheorem{prop}[thm]{Proposition}
\newtheorem{lem}[thm]{Lemma}
\newtheorem{defi}[thm]{Definition}

\newtheorem{cor}[thm]{Corollary}
\newtheorem{step}{Step}

\newtheorem{prop-defi}[thm]{Proposition-Definition}

\begin{document}
\maketitle
\begin{abstract}
In this paper, we describe the spaces of stability conditions for 
the triangulated categories associated to three dimensional 
Calabi-Yau fibrations. We deal with two cases, the 
flat elliptic fibrations and smooth $K3$ (Abelian) fibrations. 
In the first case, we will see there exist chamber structures 
similar to those of the movable cone used in birational 
geometry. In the second case, we will compare the space with the space 
of stability conditions for the closed fiber of the fibration. 
\end{abstract}

\section{Introduction}
$\quad$ In this paper, we give the descriptions of the spaces of stability 
conditions for the triangulated categories associated to
three dimensional
 Calabi-Yau fibrations. In~\cite{Brs1}, the notion of stability 
condition was introduced in order to give the mathematical framework 
for M.Douglas's work on $\Pi$-stability. For 
a Calabi-Yau variety $X$, the set of numerical stability conditions 
on $D(X)$ which satisfy local finiteness form a complex manifold, and expected 
to be an extended version of the Teichm$\ddot{\textrm{u}}$ler
 space of the stringy 
K$\ddot{\textrm{a}}$hler moduli space.
On the other hand,
the spaces of stability conditions for three dimensional 
crepant small resolutions are discussed in~\cite{Tst},
 and the resulting spaces have the
chamber structures which are similar to the chamber structure of 
Y.Kawamata's 
movable cone~\cite{Ka2}. This fact gives a connection between 
birational geometry and $\nN=2$ super conformal field theory. One  
of our purpose is to see such connections for other Calabi-Yau fibrations, 
e.g. elliptic fibrations and $K3$ (Abelian) fibrations. 

\hspace{5mm}

Let 
$S\cneq \Spec R$ for a Noetherian local complete $\mathbb{C}$-algebra $R$, 
with $0\in S$ the closed point and $\eta \in S$ the generic point.
Let
$f\colon X\to S$ be a projective surjective morphism of 
normal complex schemes with connected fibers, and $\dim X=3$, 
$\dim S \ge 1$.  
In this paper, we call it a Calabi-Yau fibration if $X$ is 
regular and $\omega _X$ is trivial. 
Note that the geometric generic fiber $X_{\bar{\eta}}$ is a Calabi-Yau 
variety in the sense that $\omega_{X_{\bar{\eta}}}$ is trivial. 
We define 
$D(X/S)$ to be the 
triangulated subcategory of $D(X)$:
$$D(X/S)\cneq \{ A\in D(X) \mid \Supp A \subset f^{-1}(0) \}.$$
Our purpose is to describe the space of stability conditions on $D(X/S)$. 
Note that $D(X/S)$ has categorical properties similar to those of 
derived categories of projective Calabi-Yau 3-folds, 
for example, the functor $E \mapsto E[3]$ gives a Serre functor.  
Thus it is interesting and important to study the stability 
conditions for $D(X/S)$. 
The case of $\dim X_{\bar{\eta}}=0$
 was partially discussed in~\cite{Brs4},~\cite{Tst},
 so we discuss the remaining cases, i.e. $f$ is an elliptic fibration 
 or a $K3$ (Abelian) fibration. 

\hspace{5mm}

At this time, we 
have to put the assumption that 
$f$ is flat if $f$ is an elliptic fibration, and 
$f$ is smooth if $f$ is a $K3$ (Abelian) fibration. 
One of the difficulties occurs when one tries 
to construct a stability condition. 
In~\cite{Brs2}, T.Bridgeland constructed 
stability conditions on the derived categories of 
$K3$ surfaces using Bogomolov inequality, 
 and it seems difficult to apply such technique without the 
 assumption above.  
 The another problem occurs, for example, if there 
 exists a projective plane $E=\mathbb{P}^2 \subsetneq f^{-1}(0)$ and 
its normal bundle is $N_{E/X}=\oO _{\mathbb{P}^2}(-3)$. 
 In this case the problem of describing stability conditions on 
 $D(X/S)$ contains the problem of describing those for 
 the total space of $\oO _{\mathbb{P}^2}(-3)$,
 which we  
 are unable to give the complete description of it
 at this time. (See~\cite{Brs4}.) 
Therefore we have to 
exclude such a divisor $E\subsetneq f^{-1}(0)$. 
 Thus we treat the cases of \textit{flat elliptic fibration} and 
  \textit{smooth K3 
 (Abelian) fibration}. 

\vspace{5mm}

For $f\colon X\to S$ as above, 
we denote by $\Stab (X/S)$
the set of locally finite numerical stability conditions on $D(X/S)$. 
Using the same argument as in~\cite{Brs2}
 and~\cite{Tst}, we can construct some standard
 points in $\Stab (X/S)$. 
Let $\beta, \omega$ be $\mathbb{R}$-divisors and assume 
$\omega$ is ample. We consider the function $Z_{(\beta, \omega)}\colon 
K(X/S) \to \mathbb{C}$ defined to be
$$Z_{(\beta, \omega)}(E)\cneq - \int e^{-(\beta +i\omega)}\ch (E)
\sqrt{\td _X}.$$
Then one can construct the t-structure with heart 
$\aA _{(\beta, \omega)} \subset D(X/S)$,  and for a suitable 
choice of $\beta, \omega$, one can check the pair 
$\sigma _{(\beta, \omega)}\cneq (Z_{(\beta, \omega)}, \aA _{(\beta, \omega)})$ 
gives a numerical stability condition on $D(X/S)$. 
We denote by $\Stab ^{\circ}(X/S)$ the connected component of $\Stab (X/S)$,
 which 
contains $\sigma _{(\beta, \omega)}$. 

\vspace{5mm}

First we give the description of $\Stab ^{\circ}(X/S)$ for flat elliptic 
fibrations. Let $V_{\mathbb{C}}
\subset N^1 (X/S)_{\mathbb{C}}$ be the $\mathbb{C}$-vector
subspace generated by $f$-vertical divisors, and 
$W_{\rm{ref}}\subset \GL (N_1 (X/S))$ be the subgroup 
generated by reflections associated to $f$-vertical divisors. 
Also let $\Lambda \subset 
\overline{\NE}(X/S)$ be the subset which consists of sums of 
extremal rational curves whose dual graphs are of Dynkin type. 
Using the 
techniques in~\cite{Tst}, we will show the following:
\begin{thm}
Let $f\colon X\to S$ be a three dimensional flat elliptic fibration, and 
assume that $f^{-1}(T)$ is smooth for a general hyper plane 
section $0\in T\subset S$. 
Then for a pair $(k,l)\in \mathbb{Z}\times N_1 (X/S)$,
one can attach a codimension two hyper plane 
$\widetilde{H}_{k,l}\subset \GL ^{+}(2, \mathbb{R})\times V_{\mathbb{C}}$ 
and has a map
$$\Stab ^{\circ}(X/S) \lr (\GL ^{+}(2, \mathbb{R}) \times V_{\mathbb{C}}) 
\setminus \bigcup _{(k,w,l)\in \mathbb{Z}\times W_{\emph{ref}} \times 
\Lambda} \widetilde{H}_{k,w(l)},$$
which is a regular covering map. 
 \end{thm}
Here it is worth recalling that $\Stab (X_{\bar{\eta}})$ is a 
universal cover of $\GL ^{+}(2, \mathbb{R})$. (See~\cite{Brs1}). 

Next we describe the 
spaces of stability conditions for $K3$ (Abelian) fibrations, 
but we will use the different kind of approach in this case. 
Our method is to compare stability conditions on $D(X/S)$ and those 
 on the derived category of the 
 special fiber $X_0 \cneq f^{-1}(0)$ studied 
in~\cite{Brs2}. According to~\cite{Brs2}, there exists a 
connected component
$\Stab ^{\circ}(X_0)$ which 
is a regular covering space over certain open subset 
$\pP _0 ^+ (X_0)\subset 
\mathbb{C} \oplus N^1 (X_0)_{\mathbb{C}} \oplus \mathbb{C}$. 
 We will show the following:
\begin{thm}
Let $f\colon X\to S$ be a three dimensional 
 smooth $K3$ fibration. Then
one has the open subset of 
$\mathbb{C}\oplus N^1 (X/S)_{\mathbb{C}}\oplus \mathbb{C}$, 
denoted by $\pP _{0}^{+}(X/S)$, and the map 
$$\zZ \colon \Stab ^{\circ}(X/S) \lr \pP _0 ^{+}(X/S)$$
which is a regular covering map. 
\end{thm}
The explicit descriptions of $\widetilde{H}_{k,l}$ and 
$\pP _0 ^{+}(X/S)$ will be given in Section 4 and Section 6 respectively. 
\subsection*{Acknowledgement}
The author thanks Tom Bridgeland for useful discussions. 
He also thanks Hokuto Uehara for reading the manuscript and 
giving the nice advice. 
He is supported by Japan Society for the Promotion of Sciences
Research Fellowships for Young Scientists, No 1611452.

\subsection*{Notations and Conventions}
For a scheme $X$, we denote 
by $\Coh (X)$ and $D (X)$ the Abelian category of 
coherent sheaves and its bounded derived category respectively. 
The shift functor on $D (X)$ is denoted by $[1]$.

\section{Stability conditions for triangulated categories}

$\quad$ In this section, we give a brief summary on   
stability conditions for triangulated categories 
introduced in~\cite{Brs1}.
We recall definitions and several properties which 
will be used in this paper. 
\subsection*{Stability conditions}
\begin{defi}
A stability condition of a triangulated category $\dD$ 
consists of a data $\sigma =(Z, \pP)$, 
where $Z\colon K(\dD)\to \mathbb{C}$ is a linear map
called a central charge, 
and full additive subcategories 
$\pP (\phi)\subset \dD$ for each $\phi \in \mathbb{R}$, 
which satisfies the following:

 \begin{itemize}
 \item $\pP (\phi +1)=\pP (\phi)[1].$ 
\item  If $\phi _1 >\phi _2$ and $A_i \in \pP (\phi _i)$, then 
$\Hom (A_1, A_2)=0$. 
\item  If $0\neq E\in \pP (\phi)$, then $Z(E)=m(E)\exp (i\pi \phi)$ for some 
$m(E)\in \mathbb{R}_{>0}$. 
\item (Harder-Narasimhan filtration)
For a non-zero object $E\in \dD$, we have the 
following collection of triangles:
$$\xymatrix{
0=E_0 \ar[rr]  & &E_1 \ar[dl] \ar[rr] & & E_2 \ar[r]\ar[dl] & \cdots \ar[rr] & & E_n =E \ar[dl]\\
&  A_1 \ar[ul]^{[1]} & & A_2 \ar[ul]^{[1]}& & & A_n \ar[ul]^{[1]}&
}$$
such that $A_j \in \pP (\phi _j)$ with $\phi _1 > \phi _2 > \cdots >\phi _n$. 
\end{itemize}
\end{defi} 
We can see each $\pP (\phi)$ is an Abelian category,
and the non-zero objects of $\pP (\phi)$ are called semistable of 
phase $\phi$. 
The objects $A_j$ are called semistable factors of $E$ with 
respect to $\sigma$, and 
we write $\phi _{\sigma}^{+}(E)=\phi _1$ and $\phi _{\sigma}^{-}(E)=\phi _n$. 
It is an easy exercise to check that the decompositions into 
semistable factors $A_i$ are unique up to isomorphism. 
In particular if there exists another stability condition $\sigma '=(Z', \pP')$
with $\pP (\phi)\subset \pP '(\phi)$ for all $\phi \in \mathbb{R}$, 
then $\pP (\phi)=\pP '(\phi)$. 
In this paper, we introduce the notation $\pP _s (\phi)$ to be 
$$\pP _s (\phi)\cneq \{ E\in \pP (\phi) \mid E\mbox{ is a simple object of }
\pP (\phi)\}.$$
The objects of $\pP _s (\phi)$ are called stable. 
The mass of $E$ is defined to be 
$$m_{\sigma}(E)=\sum _{j}|Z(A_j)|.$$
For an interval $I\subset \mathbb{R}$, denote by $\pP (I)$ the minimum
extension closed subcategory of $\dD$ which contains 
 $\pP (\phi)$ for $\phi \in I$. 
If $I=(a,b)$ with $b-a \le 1$, then $\pP ((a,b))$ is a quasi-Abelian 
category and $\pP ((0,1])$ is an 
Abelian category. (In quasi-Abelian category, one has kernel and cokernel, 
but image and coimage may not coincide.)
In fact, 
$\pP ((0,1])$ is a heart of some t-structure on $\dD$. 
This construction provides the following proposition:
\begin{prop}\emph{\bf{\cite[Proposition 4.2]{Brs1}}}\label{tst}
Giving a stability condition on $\dD$ is equivalent to giving a 
bounded t-structure on $\dD$ with heart $\aA$, 
and a linear function
$Z\colon K(\dD)\to \mathbb{C}$
such that 
$$0\neq E \in \aA \ \Rightarrow \ Z(E)\in \mathbb{R}_{>0}\exp (i\pi \phi)
\mbox{ with } 0<\phi \le 1,$$
and the pair $(Z, \aA)$ satisfies Harder-Narasimhan property. 
\end{prop}

We have to 
put the locally finiteness conditions to introduce 
the topology on the set of stability conditions.  This means 
for each $\phi \in \mathbb{R}$, there exists $\varepsilon >0$ such 
that quasi-Abelian category 
$\pP ((\phi -\varepsilon, \phi +\varepsilon))$ is of finite length, i.e. 
noetherian and artinian with respect to the strict monomorphism.
 (\textit{strict} means image and coimage coincide.)
In particular each $\pP (\phi)$ is also of finite length, hence
has a Jordan-H\"{o}rder decomposition. 
We denote by $\Stab (\dD)$ the set of locally finite  
stability conditions on $\dD$. 
Forgetting the information of $\pP$, we have the map:
$$\zZ \colon \Stab (\dD) \longrightarrow \Hom _{\mathbb{Z}}(K(\dD), 
\mathbb{C}).$$
We can induce the natural topology on $\Stab (\dD)$ 
so that the map $\zZ$ becomes continuous. 
If $\sigma$ and $\tau$ are stability conditions on a triangulated 
category $\dD$, then define $d(\sigma, \tau)$ to be 
$$d(\sigma, \tau)\cneq \sup _{E\neq 0}\left\{ \lvert \phi _{\tau}^{-}(E)
-\phi _{\sigma}^{-}(E)\rvert, \lvert \phi _{\tau}^{+}(E)
-\phi _{\sigma}^{+}(E)\rvert \right\} \in [0, \infty]. $$
Also for $\sigma =(Z, \pP) \in \Stab (\dD)$, we can induce the generalized 
norm $\lVert \cdot \rVert _{\sigma}$
on $\Hom _{\mathbb{Z}}(K(\dD), \mathbb{C})$,
\begin{align*}\lVert U \rVert _{\sigma} &\cneq \sup \left\{
\frac{\lvert U(E) \rvert}{\lvert Z(E) \rvert} :
E\mbox{ is semistable in }\sigma \right\} \\
& = \sup \left\{
\frac{\lvert U(E) \rvert}{\lvert Z(E) \rvert} :
E\mbox{ is stable in }\sigma \right\}.\end{align*}
Then the following subsets of $\Stab (\dD)$, 
 $$B_{\varepsilon}(\sigma) \cneq \{
 \tau \in \Stab (\dD) \mid d(\sigma, \tau)<\varepsilon, 
 \lVert \zZ (\tau)-\zZ (\sigma) \rVert _{\sigma} < \sin \pi \varepsilon \}$$
 provides an open basis of $\Stab (\dD)$. 
 One can see the norms $\lVert \cdot \rVert _{\sigma}$ and $\lVert \cdot \rVert _{\tau}$
are equivalent if $\sigma$, $\tau$ are contained in the same connected 
component of $\Stab (\dD)$. Thus for each connected component 
$\Sigma \subset \Stab (\dD)$, we have the well-defined topology on 
$\Hom _{\mathbb{Z}}(K(\dD), \mathbb{C})$. 
Let $V(\Sigma)\subset \Hom _{\mathbb{Z}}(K(\dD), \mathbb{C})$
be
$$V(\Sigma)\cneq \{ U\in \Hom _{\mathbb{Z}}(K(\dD), \mathbb{C}) :
\lVert U \rVert _{\sigma} <\infty \mbox{ for } \sigma \in \Sigma \}.$$
\begin{thm}\emph{\bf{\cite[Theorem 1.2]{Brs1}}} \label{lhom}
For each connected component $\Sigma \subset \Stab (\dD)$, 
$\zZ$ restricts to give a local homeomorphism, 
$\zZ \colon \Sigma \to V(\Sigma)$. 
\end{thm}

One of the key lemma for the proof of Theorem~\ref{lhom} is the following
deformation result,  
which we will use in Section 5. 
\begin{thm}\emph{\bf{\cite[Theorem 6.1]{Brs1}}} \label{def}
Let $\sigma =(Z, \pP)$ be a locally finite stability condition 
on $\dD$. Let us take $0<\epsilon _0 \le 1/8$ such that 
$\pP ((t-4\epsilon _0, t+4\epsilon _0))\subset 
\dD$ is of finite length for all $t\in \mathbb{R}$. 
Then if 
$0<\epsilon  <\epsilon _0$ and $W\colon K(\tT)\to \mathbb{C}$ is 
a linear map satisfying
$$ \lVert W-Z \rVert _{\sigma} <\sin (\pi \epsilon), $$
then there exists a stability condition $\tau =(W, \qQ)$ on $\dD$ 
with $d(\sigma, \tau)<\epsilon$. 
In particular if $\Imm \zZ \subset \mathbb{C}$ is a discrete
subgroup, one can take $\varepsilon _0$ to be $1/8$. 
\end{thm}

We have the action of the group of autoequivalences $\Aut (\dD)$
to $\Stab (\dD)$: for $\Phi \in \Aut (\dD)$ and $\sigma =(Z,\pP)$, 
$\Phi (\sigma)=(Z', \pP ')$ with $\pP ' (\phi)=\Phi (\pP (\phi))$ 
and $Z'(E)=Z(\Phi ^{-1}(E))$. 
Also the additive group $\mathbb{C}$ acts on $\Stab (\dD)$: 
for $\lambda \in \mathbb{C}$, 
$\lambda (\sigma)=(Z'',\pP '')$ with $\pP '' (\phi )=\pP (\phi +\Ree \lambda)$
and $Z'' (E)=\exp (-i\pi \lambda)Z(E)$. This action commutes with the 
action of autoequivalences.

\subsection*{Wall and chamber structures}

We recall some facts discussed in~\cite{Brs2}
on the existence of wall and chamber 
structures on the space of stability conditions. 
In general $\Stab (\dD)$ may be infinite dimensional. Therefore 
in usual we consider stability conditions $\sigma =(Z, \pP)$ 
such that $Z$ factors through the surjection $K(\dD)\twoheadrightarrow \nN$
for a fixed finitely generated $\mathbb{Z}$-module $\nN$.
(See numerical stability conditions in~\cite{Brs2}.)
Let $\Stab _{\nN}(\dD)$ be the set of locally finite stability 
conditions $\sigma =(Z, \pP)$ such that $Z$ factors $\nN$, 
$Z\colon K(\dD) \to \nN \to \mathbb{C}$. 
Then each connected component $\Sigma \subset \Stab _{\nN}(\dD)$ 
carries a map into $\nN _{\mathbb{C}}^{\ast}$, thus Theorem~\ref{lhom}
implies $\Sigma$ is a complex manifold. 
For each $m\in \mathbb{N}$
and $\sigma \in \Sigma$, we denote by $\sS _m$ the set of 
objects, 
$$\sS _m \cneq \{ E\in \dD \mid m_{\sigma}(E) <m \}.$$
Then let us consider the following condition $(\diamondsuit)$, 
$$(\diamondsuit) \quad 
\mbox{ For each }m\in \mathbb{N}, 
\mbox{ the set } \{ [E] \in \nN \mid E\in \sS _m \}
\mbox{ is finite.}$$
Note that if $(\diamondsuit)$ holds for some $\sigma$, 
then it holds for 
every points in $\Sigma$. 
The following proposition is due to~\cite[Proposition 8.3]{Brs2}.
\begin{prop}\label{wall}
Assume the condition $(\diamondsuit)$ holds for $\sigma \in \Sigma$, and 
let $\sS$ be the subset of $\sS _m$ for some $m$. Then
for a
fixed compact subset $O\subset \Sigma$, there is a finite 
number of real codimension one submanifolds 
$\{ \wW_{\gamma} \mid \gamma \in \Gamma \}$ 
such that each connected component  
$$O^{\circ}\subset O\setminus \bigcup _{\gamma \in \Gamma}\wW_{\gamma}$$
has the following property: If $E\in \sS$ is semistable in $\sigma$ for 
some $\Sigma$,
 then $E$ is semistable in $\sigma$ for all $\sigma 
\in O^{\circ}$. If $[E]\in \nN$
 is primitive, then $E$ is in fact stable.
 \end{prop}
\textit{Proof}. 
The same proof of~\cite[Proposition 8.3, Corollary 8.4]{Brs2} is applied.

\subsection*{Numerical stability conditions for Calabi-Yau fibrations}
Let $f\colon X\to S$ be a Calabi-Yau fibration as in Introduction,  
$X_0 \cneq f^{-1}(0)$ and 
denote by $D(X/S)$ to be the subcategory of $D(X)$, 
$$D(X/S)\cneq \{ E\in D(X) \mid \Supp (E)\subset X_0 \}.$$
Then $D(X/S)$ is an Ext-finite category, has a Serre functor $S_X=[3]$. 
Let $K(X/S)$ be the Grothendieck group of $D(X/S)$. We have the 
following pairing:
$$\chi \colon 
K(X/S)\times K(X) \ni (E, F) \longmapsto \chi (E, F)\cneq \sum (-1)^i \dim 
\Ext ^i (E, F) \in \mathbb{Z}.$$
We say $E_1, E_2 \in K(X/S)$ (resp $E_1, E_2 \in K(X)$)
 are numerically equivalent if 
$\chi (E_1, F)=\chi (E_2, F)$ for all $F \in K(X)$ (resp $ F\in K(X/S)$). 
Then define $\nN (X/S)$ and $\nN (X)$ to be the numerical equivalence 
classes of $K(X/S)$, $K(X)$ respectively. Thus $\chi$ descends to 
the perfect pairing,
$$\chi \colon \nN (X/S) 
\times \nN (X) \lr \mathbb{Z}.$$
For a stability condition $\sigma =(Z, \pP)$ on $D(X/S)$, 
we call it \textit{numerical} if $Z$ factors though the 
surjection $K(X/S) \twoheadrightarrow \nN (X/S)$. We denote by $\Stab (X/S)$
the set of locally finite numerical stability conditions. By restricting 
$\zZ$, we obtain the map, 
$$\zZ \colon \Stab (X/S) \lr \Hom _{\mathbb{Z}}
(\nN (X/S), \mathbb{C})\cong \nN (X)_{\mathbb{C}}.$$
Let us describe $\nN (X)_{\mathbb{C}}$ explicitly. 
We denote by $\CH ^p (X)$ the Chow group, which 
 is a rational equivalence class of codimension $p$-cycles on $X$. 
 Then we have the following map,
 $$v\colon K(X) \ni E \longmapsto \ch (E)\sqrt{\td _X} \in
  \bigoplus _{p\ge 0} \CH ^p(X)_{\mathbb{Q}}, $$
 which is taking Mukai vectors. On the other hand,  
we define the $\mathbb{R}$-vector 
 space $N^p (X/Y)$ to be the numerical equivalence
 classes of 
 codimension $p$ cycles:
 $$N^p (X/Y) \cneq \bigoplus _{D\subset X} \mathbb{R}[D] / \equiv .
 $$
 Here for codimension $p$-cycles $D_1, D_2$ on $X$, 
  $D_1 \equiv D_2$ if and only 
 if $D_1 \cdot Z=D_2 \cdot Z$ for any $p$-dimensional cycle $Z\subset 
 X_0$.
 Then let $\pi ^p \colon \CH ^p (X) \to N^p (X/Y)$
 be the natural map when $p\le \dim X_0$ and zero map 
 when $p>\dim X_0$. Then the composition 
 $$\pi \circ v \colon K(X) \stackrel{v}{\lr} \bigoplus _{p\ge 0}
 \CH (X) _{\mathbb{Q}} \stackrel{\oplus \pi ^p}{\lr} 
 \bigoplus _{p=0}^{\dim X_0}N^p (X/Y),$$
factors though the quotient $K(X)\twoheadrightarrow \nN (X)$. 
In fact if $E \in K(X)$ is numerically trivial, then 
$\ch _p (E) =0$ for $p \le \dim X_0$. Thus
we have an isomorphism, 
 $$\nN (X)_{\mathbb{C}} \stackrel{\cong}{\lr}
 \bigoplus _{p=0}^{\dim X_0} N^p (X/S)_{\mathbb{C}}.$$
If $f$ is a flat elliptic fibration or smooth $K3$ (Abelian) fibration, 
we have $N^i(X/S)=\mathbb{C}$ for $i=0, \dim X_0$. 
Therefore we have the following
maps:
$$\left\{
\begin{array}{lll}
\zZ \colon \Stab (X/S)  \lr \mathbb{C} \oplus N^1 (X/S)_{\mathbb{C}}
& (\dim X_{0}=1) \\
\zZ \colon \Stab (X/S) \lr \mathbb{C} \oplus N^1 (X/S)_{\mathbb{C}}
\oplus \mathbb{C} & (\dim X_0 =2)
\end{array}\right.$$
It is well-known that $N^1 (X/S)$ is finite dimensional. So each connected 
component of $\Stab (X/S)$ is a complex manifold.

\section{Calabi-Yau fibrations and cone structures}
Here we recall some terminologies from birational 
geometry, especially used in~\cite{Ka2}.
Let $f\colon X\to S$ be a Calabi-Yau fibration with 
$0 \in S$ the closed point and $\bar{\eta}\to S$ the 
geometric generic point. 
We use the following Cartesian squares, 
$$\xymatrix{
X_{\bar{\eta}} \ar[r] ^{j} \ar[d] & X \ar[d] & X_0 \ar[l] _{i} \ar[d] \\
\bar{\eta} \ar[r] & S & \ar[l] 0. }$$
Then a Cartier divisor $D$ on $X$ is called 
\begin{itemize}  
\item $f$-\textit{big} if the Kodaira dimension 
$\kappa (X_{\bar{\eta}}, j^{\ast}D)$ is equal to 
$\dim X_{\bar{\eta}}$.
We denote by 
$\bB (X/S)\subset N^1 (X/S)$ the open convex cone 
generated by $f$-big divisors.
\item $f$-\textit{nef} if $D\cdot C \ge 0$ for any curve $C$ on $X_0$. 
The nef cone $\overline{\aA}(X/S)\subset N^1 (X/S)$ is the closed convex cone 
generated by $f$-nef divisors. The set of 
its inner points are denoted by $\aA (X/S)\subset 
\overline{\aA}(X/S)$, which consists of numerical classes of 
$\mathbb{R}$-ample divisors.  
\item $f$-\textit{movable} if 
$$\codim \Supp  (f^{\ast}f_\ast \oO _X (D) \lr \oO _X (D))\ge 2.$$
The movable cone 
$\overline{\mM}(X/S)\subset N^1 (X/S)$ is a closed convex cone 
generated by $f$-movable divisors.
\item $f$-\textit{vertical} if $f(D)\subsetneq S$. We denote by 
$V\subset N^1 (X/S)$ the vector subspace generated by 
$f$-vertical divisors.
\end{itemize}
We have the following inclusions:
$$\overline{\aA}(X/S)\subset \overline{\mM}(X/S)\subset
\overline{\bB}(X/S)\subset N^1 (X/S).$$
On the other hand, let $N_1 (X/S)$ be the $\mathbb{R}$-vector 
space generated by numerical classes of 
one cycles on $X_0$. 
Note that we have the perfect pairing, 
$$N^1 (X/S) \times N_1 (X/S) \ni (D, C) \longmapsto 
D\cdot C \in \mathbb{R}.$$
Let
$\overline{\NE}(X/S)\subset N_1 (X/S)$ be
the closed cone generated by effective one cycles. Thus 
$\overline{\NE} (X/S)$ is a dual cone of $\overline{\aA}(X/S)$ under the 
above pairing.

For another Calabi-Yau fibration $f'\colon W\to S$ with a $S$-birational 
map $\phi \colon W\dashrightarrow X$,
let $\phi _{\ast} \colon N^1 (W/S) \to N^1 (X/S)$ be the
strict transform. 
We have the following lemma from~\cite{Ka2}.
\begin{lem}\emph{\bf{\cite{Ka2}}} \label{nef}
For two Calabi-Yau fibrations 
$f_i '\colon W_i \to S$ with $S$-birational maps 
$\phi _i \colon W_i \dashrightarrow X (i=1,2)$, 
if $\phi _{1\ast}\aA (W_1/S)\cap \phi _{2\ast}\aA (W_2/S)\neq \emptyset$, 
then $\phi _1 ^{-1}\circ \phi _2 
\colon W_2 \dashrightarrow W_1$ is an isomorphism. 
\end{lem}
In case of there dimensional flat elliptic fibration, we have the following 
lemma on the structures of those cones.  
\begin{lem}\emph{\bf{\cite{Ka2}}}
\label{mov}
Let $f\colon X\to S$ be a three dimensional flat elliptic fibration. Then 

(1) The nef cone 
$\overline{\aA}(X/S)$ is a rational polyhedral cone
in $N^1 (X/S)$. 

(2)  We have 
the following decomposition
$$\overline{\mM}(X/S)\cap \bB (X/S) =\bigcup _{\phi \colon W\dashrightarrow X}
\phi _{\ast}
\overline{\aA} (W/S) \cap \bB (X/S),$$
which is locally finite inside $\bB (X/S)$. 

(3) The open cone  $\bB (X/S)$ is generated by $\overline{\mM}(X/S)\cap 
\bB (X/S)$ and $f$-vertical divisors. 
\end{lem}

\section{Stability conditions for flat elliptic fibrations}
In this section, we assume $f\colon X\to S$ is a 
three dimensional flat elliptic fibration, and 
give the description of the spaces of stability 
conditions on $D(X/S)$. 
The strategy is almost same as in~\cite{Tst}, and we will leave some 
detailed discussions to the readers. For a cone $\aA \subset N^1 (X/S)$, 
we denote by $\aA _{\mathbb{C}}$ its complexified cone, 
$$\aA _{\mathbb{C}}\cneq \{ \beta +i\omega \in N^1 (X/S)_{\mathbb{C}} \mid
\omega \in \aA\}.$$ 

\subsection*{The construction of stability conditions}
For $\beta +i\omega \in N^1 (X/S)_{\mathbb{C}}$, let 
$Z_{(\beta, \omega)}\colon \nN (X/S) \to \mathbb{C}$ be
\begin{align*}
Z_{(\beta, \omega)}(E) & \cneq -\int e^{-(\beta +i\omega)}\ch (E) 
\sqrt{\td _X} \\
& = - \ch _3 (E) +(\beta +i\omega)\ch _2 (E). 
\end{align*}
Under the isomorphism $\Hom (\nN (X/S), \mathbb{C})\cong \mathbb{C} 
\oplus N^1 (X/S)_{\mathbb{C}}$, $Z_{(\beta, \omega)}$ corresponds to 
the element $(-1, \beta +i\omega)$. 
Let $\Coh (X/S)\cneq \Coh (X)\cap D(X/S)$. As in~\cite[Lemma 4.1]{Tst}, 
 we have the following 
lemma:
\begin{lem}\label{vol}
For $\beta +i\omega \in \aA (X/S)_{\mathbb{C}}$,
the pair $\sigma _{(\beta, \omega)}\cneq (Z_{(\beta, \omega)}, \Coh (X/S))$
determines a point of $\Stab (X/S)$.
\end{lem}
\textit{Proof}. We can apply the same proof of Lemma~\cite[Lemma 4.1]{Tst}. 
$\quad \square$
 
 \vspace{5mm}

It is easy to see that the stability conditions $\sigma _{(\beta, \omega)}$ 
constructed as above are contained in the same connected component 
denoted by 
$\Stab ^{\circ}(X/S) \subset \Stab (X/S)$.  
We have the following:
\begin{lem}\label{di}
For $\sigma _{(0, \omega)}\in \Stab (X/S)$,
 the condition $(\diamondsuit)$ in Section 2 
 is satisfied. 
Thus there exist wall and chamber structure on the connected 
component $\Stab ^{\circ}(X/S)$. 
\end{lem}
\textit{Proof}. 
Fix $m>0$ and consider the set 
$\sS _m \cneq \{ E\in D(X/S) \mid m_{\sigma _{(0, \omega)}}(E) <m \}$.
 Then for $E\in \sS _m$, 
 the stable factors of $E$ in $\sigma _{(0, \omega)}$ 
 are also contained in $\sS _m$. Since stable 
 objects in $\sigma _{(0, \omega)}$ are shift of stable sheaves, 
 it suffices to show the numerical classes of $\sS _m \cap \Coh (X/S)$ 
 are finite.
 For $E\in \sS _m \cap \Coh (X/S)$, 
 note that the numerical class of $E$ is determined by the 
 pair $(\ch _2 (E), \ch _3 (E))\in N_1 (X/S)\oplus \mathbb{R}$. 
 Let
$C \subset X_0$ be one of the irreducible 
components of $X_0$. Since $|\omega \cdot \ch _2 (E)|<m$, 
the generic length of $E$ at $C$ is not greater than $m/(C \cdot \omega)$. 
Also since $|\ch _3 (E)|<m$, we can conclude 
there exists finite number of possibilities for the pair $(\ch _2 (E), 
\ch _3 (E))$.   $\quad \square$.

\subsection*{Normalized stability conditions}
Let us define $\Stabn (X/S)$ to be the normalization 
of $\Stab (X/S)$ under the action of $\mathbb{C}$,
$$\Stabn (X/S) \cneq \{ \sigma =(Z, \pP)\in \Stab (X/S) 
\mid Z([\oO _x])=-1 \}.$$
Note that for $\sigma =(Z, \pP)\in \Stabn (X/S)$, $Z$ is written 
as $Z_{(\beta, \omega)}$ for some $\beta +i\omega \in N^1 (X/S)_{\mathbb{C}}$.
Restricting $\zZ$ to $\Stabn (X/S)$, we obtain the
map,
$$\zZ_{\rm{n}} \colon 
\Stabn (X/S) \lr N^1 (X/S)_{\mathbb{C}},$$
such that $\zZ _{\rm{n}}$ takes $\sigma =(Z_{(\beta, \omega)}, \pP)$
to $\beta +i\omega$.  
Note that 
all the stability conditions 
$\sigma _{(\beta, \omega)}$ constructed in Lemma~\ref{vol} are 
contained in $\Stabn (X/S)$, and define $U_X$ to be 
$$U_X \cneq \{ \sigma _{(\beta, \omega)} \in \Stabn (X/S) 
\mid \beta +i\omega \in 
\aA (X/S)_{\mathbb{C}} \}.$$ 
Then $U_{X}$ is an open subset of $\Stabn (X/S)$,  and
$\zZ_{\rm{n}}$ restricts to give a homeomorphism 
between $U_X$ and $\aA (X/S)_{\mathbb{C}}$. 
Let $\Stabn ^{\circ}(X/S)$ be the connected component 
of $\Stabn (X/S)$ which 
contains $U_X$. Note that we have the inclusion
$$\Stabn ^{\circ}(X/S) \subset \Stab ^{\circ}(X/S) \cap \Stabn (X/S).$$

\subsection*{Other regions by Fourier-Mukai transform}
Let $f'\colon W\to S$ be another three dimensional 
flat elliptic fibration. We say the equivalence 
$\Phi \colon D(W)\to D(X)$
is of birational 
Fourier-Mukai type over $S$ if there exists an object 
$\rR \in D(X\times W)$, which is supported on $X\times _S W$ such 
that $\Phi$ is written as 
$$\Phi = \Phi _{W\to X}^{\rR}\cneq 
\dR p_{X\ast}(p_W ^{\ast}(\ast)\dotimes \rR) \colon D(W)\to D(X), $$
and if we consider the associated functor between derived 
categories of quasi-coherent sheaves
$\Phi _{W\to X}^{\rR} \colon D(\QCoh (W)) \lr D(\QCoh (X))$,
then it takes $\oO _{k(W)}$ to $\oO _{k(X)}$. 
Here $p_X$, $p_W$ are projections from $X\times W$ onto corresponding 
factors, and $k(X)$, $k(W)$ are generic points of $X$, $W$ respectively.  
Note that $\Phi$ induces the $S$-birational map 
$\phi \colon W \dashrightarrow X$. 
\begin{defi}
We define the set $\FM _{\rm{bir}} (X)$ to be the equivalence class of 
pairs $(W, \Phi)$ such that $g\colon W\to S$ is another flat elliptic
 fibration and $\Phi \colon 
D(W)\to D(X)$ is of birational 
 Fourier-Mukai type over $S$. 
\end{defi}

For $(W, \Phi)\in \FM _{\rm{bir}}(X)$, let 
$\phi \colon W\dashrightarrow X$ be the induced $S$-birational map and 
define $\widetilde{\phi}$ to be the map
$$\widetilde{\phi}\colon N^1 (W/S)_{\mathbb{C}} \ni 
\beta +i\omega \longmapsto c_1 (\Phi (\oO _W))+\phi _{\ast}\beta 
+i\phi _{\ast}\omega \in N^1 (X/S).$$
The following proposition is a direct adaptation 
of~\cite[Proposition 4.8]{Tst} in our case.  
\begin{prop}\label{com}
Take $(W, \Phi)\in \FM _{\emph{bir}} (X)$. Then $\Phi$ induces a 
homeomorphism $\widetilde{\Phi}\colon \Stabn (W/S)\to \Stabn (X/S)$ 
which fits into the commutative diagram, 
$$\xymatrix{
\Stabn (W/S) \ar[r]^{\widetilde{\Phi}}\ar[d]_{\zZ _{\rm{n}}} & \Stabn (X/S) \ar[d]^{\zZ _{\rm{n}}} \\
N^1 (W/S)_{\mathbb{C}} \ar[r]^{\widetilde{\phi}}& N^1 (X/S)_{\mathbb{C}}.
}$$
\end{prop}
\textit{Proof}. The proof is same as in~\cite[Proposition 4.8]{Tst}. 
$\quad \square$

\vspace{5mm}

Now for $(W, \Phi)\in \FM _{\rm{bir}}(X)$,
 we construct the region $U(W, \Phi)$ 
to be 
$$U(W, \Phi)\cneq \widetilde{\Phi}(U_W) \subset \Stabn (X/S).$$
By proposition~\ref{com}, $\zZ _{\rm{n}}$ restricts to give a homeomorphism 
between $U(W, \Phi)$ and $\phi _{\ast}\aA (W/S)_{\mathbb{C}}$.

\subsection*{Codimension one boundaries of $U_X$}
We study the codimension one boundaries of $U_X$
in $\Stabn ^{\circ}(X/S)$. 
Let us take $\sigma =(Z_{(\beta, \omega)}, \pP)\in \overline{U}_X$. 
Then $\beta +i\omega$ lies in the nef cone $\overline{\aA}(X/S)_{\mathbb{C}}$. 
Note that since $f$ is flat, $\overline{\aA}(X/S)$ is a rational 
polyhedral cone by Lemma~\ref{mov}. We say $\sigma$ lies in the codimension 
one boundary if and only if $\beta +i\omega$ lies in the codimension one 
wall of $\overline{\aA}(X/S)_{\mathbb{C}}$. 

\begin{lem}\label{bound}
For $\sigma =(Z_{(\beta, \omega)}, \pP)
\in \overline{U}_X$, we have 
$\beta +i\omega 
\in \overline{\aA}(X/S)_{\mathbb{C}}\cap \bB (X/S)_{\mathbb{C}}$. 
\end{lem}
\textit{Proof}.
First consider the morphism $j\colon X_{\bar{\eta}} \to X$ and 
the pull-back $j^{\ast}\colon K(X) \to K(X_{\bar{\eta}})$. 
If $E\in K(X)$ is numerically zero, then $\ch _p (E)=0$ for 
$p=0, 1$. Therefore $\ch _p (j^{\ast}E)=0$ for $p=0,1$ and this 
implies $j^{\ast}E \in K(X_{\bar{\eta}})$ is also numerically zero. 
Thus $j^{\ast}\colon K(X)\to K(X_{\bar{\eta}})$ descends to the map
$j^{\ast} \colon\nN (X) \to \nN (X_{\bar{\eta}})$, 
and we obtain the map $j_{\ast}\colon \nN (X_{\bar{\eta}})
\to \nN (X/S)$ by taking the dual. Note that we can identify 
$\nN (X_{\bar{\eta}})$ and $\mathbb{Z}^{\oplus 2}$ by the 
map, 
$$\nN (X_{\bar{\eta}}) \ni E \longmapsto (\rk E, \deg E)\in \mathbb{Z}\oplus 
\mathbb{Z}.$$
Take $(r, d)\in \mathbb{Z}^{\oplus 2}\in \nN (X_{\bar{\eta}})$ such 
that $r$ and $d$ are coprime, and $r>0$. Then for 
$\sigma '=(Z_{(\beta ', \omega ')}, \Coh (X/S))\in U_X$ with 
$\beta '+i\omega ' \in \aA (X/S)_{\mathbb{Q}}$, 
we can consider the relative moduli theory of $(\beta', \omega ')$-twisted
semistable sheaves 
$E \in\Coh (X/S)$ 
with $[E] =j_{\ast}(r,d) \in \nN (X/S)$. 
Let
$\overline{\mM}(r,d)\to S$ be its coarse moduli space. 
Its geometric generic fiber $\overline{\mM}(r,d)_{\bar{\eta}}$
 is nothing but the 
moduli space of $j^{\ast}(\beta ', \omega ')$-twisted 
semistable sheaves on $X_{\bar{\eta}}$, which is non-empty. Since
$\overline{\mM}(r,d)$ is projective over $S$, it follows that 
the closed fiber $\overline{\mM} _0 (r,d)$ is also non-empty. Thus for 
each $\sigma '
\in U_X$, there exists $E\in \Coh (X/S)$ with 
$[E]=j_{\ast}(r,d)\in \nN (X/S)$ such that $E$ is semistable in
$\sigma '$. By Lemma~\ref{di},
for each 
$\sigma \in \partial{U}_X$
there exists $E\in D(X/S)$ which is semistable in $\sigma$ and 
$[E]=j_{\ast}(r,d)\in \nN (X/S)$.  Therefore the composition 
$$ 
Z_{(\beta, \omega)} \circ j_{\ast} \colon \nN (X_{\bar{\eta}}) \lr 
\nN (X/S) \lr \mathbb{C}$$
does not have kernel. Let us assume $\deg (\omega |_{X_{\bar{\eta}}})=0$. 
Then $\omega =0$ since $\omega$ is nef. 
We may also assume $\beta$ is rational, hence $\beta |_{X_{\bar{\eta}}}=d/r 
\in N^1 (X_{\bar{\eta}})_{\mathbb{Q}}$ for some $(r,d)$ which are coprime
and $r>0$.
Therefore if we take 
 $E\in \nN (X_{\bar{\eta}})$
 with $(\rk E, \deg E)=(r,d)$, we have 
\begin{align*}
Z_{(\beta, \omega)}(j_{\ast}[E]) &= Z_{(\beta |_{X_{\bar{\eta}}}, 
\omega |_{X_{\bar{\eta}}})}(E) \\
&= -d +(\beta |_{X_{\bar{\eta}}})\cdot r \\
& =0, \end{align*}
which is a contradiction. $\quad \square$

\vspace{5mm}

Now we have proved $\zZ _{\rm{n}}(\overline{U}_X)\subset 
\overline{\aA}(X/S)_{\mathbb{C}}\cap \bB (X/S)_{\mathbb{C}}$.
Any element $\omega \in \overline{\aA}(X/S)\cap \bB (X/S)$ corresponds to 
the birational contraction $g\colon X\to Y$
 with $\omega =g^{\ast}\omega _Y$ for $\omega _Y \in 
 \aA (Y/S)$. Note that the dimension of any fiber of $g$ is
less than or equal to one. In this case, one can construct the 
heart of perverse t-structure 
$^{d} \PPer (X/Y)$ in the sense of~\cite{Br1}. 
To introduce this, define $\cC\cneq \{ c\in \Coh (X) \mid \dR g_{\ast}c=0 \}$
and 
$$^{d} \PPer (X/Y)\cneq 
\{ E\in D(X) \mid \dR g_{\ast}E\in \Coh (Y), \Hom ^p(\cC, E)
=\Hom ^p(E, \cC)=0
, p<-d \}.$$
Assume $\omega \in \overline{\aA}(X/S)\cap \bB (X/S)$ is contained in 
the codimension one wall and $l \subset \overline{\NE}(X/S)$ is 
the extremal ray with supporting function $\omega$. 
We have the following two types:
\begin{itemize}
\item Type I: $\omega \in \overline{\aA}(X/S)\cap \bB (X/S)$ is in 
Type I wall if and only if there exists a diagram 
$$\xymatrix{
(C\subset X) \ar[r]^{g} \ar[rd]_{f} & (0\in Y) \ar[d]^{h} & \ar[l]_{g^{\dag}}
 (X^{\dag}\supset 
C^{\dag}) \ar[ld]^{f^{\dag}} \\
& S &, }$$
such that $\omega =g^{\ast}\omega _Y$
with $\omega _Y$ ample on $Y$. 
Here $g$ is a flopping contraction which contracts 
only single rational curve $C$ and $X^{\dag}\dashrightarrow X$ is its flop.
In this case, we have the equivalence~\cite{Br1},~\cite{Ch},
$$\Phi _l \cneq \Phi _{X^{\dag} \to X}^{\oO _{X\times _Y X^{\dag}}}
\colon D(X^{\dag}) \lr D(X),$$
which takes $\iPPer (X^{\dag}/Y)$ to $\oPPer (X/Y)$.  
In~\cite{Tst}, we called such equivalence as \textit{standard}, 
and the corresponding isomorphism $\widetilde{\phi}_l \colon 
N^1 (X^{\dag}/S)_{\mathbb{C}} \to N^1 (X/S)_{\mathbb{C}}$ is the strict 
transform for the birational map $g^{-1}\circ g^{\dag} \colon X^{\dag} 
\dashrightarrow X$. 

\item Type II: $\omega \in \overline{\aA}(X/S)\cap \bB (X/S)$ is in 
Type II wall if and only if there exists a diagram 
$$\xymatrix{
(E\subset X) \ar[r]^{g} \ar[rd]_{f} & (Z\subset Y) \ar[d]^{h} \\
& S, }$$
such that $g$ is a divisorial contraction whose restriction to $E$ 
is a $\mathbb{P}^1$-bundle, $g|_{E}\colon E\to Z$, and 
$\omega =g^{\ast}\omega _Y$
with $\omega _Y$ ample on $Y$.  
In this case the moduli space of perverse point sheaves in 
the sense of~\cite{Br1} is $X$ itself, and one has the
autoequivalence 
$$\Phi _l \colon D(X) \lr D(X),$$
which takes $\iPPer (X/Y)$ to $\oPPer (X/Y)$. 
$\Phi _l$ is written as an $EZ$-spherical twist introduced in~\cite{Ho}. 
The corresponding isomorphism $\widetilde{\phi}_l \colon N^1 (X/S)_{\mathbb{C}}
\to N^1 (X/S)_{\mathbb{C}}$ is written as the reflection
$$\widetilde{\phi}_l (\beta)=\beta +(\beta \cdot l)[E].$$
\end{itemize}

In both cases, 
let $C\subset X_0$ be an irreducible rational curve 
which generates an extremal ray $l\subset \overline{\NE}(X/S)$.
Let us
take $\lL \in \Pic (X)$ such that $\lL \cdot C=1$. 
As in~\cite{Tst}, we have the following lemma:
\begin{lem}\label{bou}

Assume $\beta +i\omega \in \overline{\aA}
 (X/S)_{\mathbb{C}}\cap \bB (X/S)_{\mathbb{C}}$
lies in the codimension one wall, and 
let $g\colon X \to Y$ be the corresponding birational contraction. 
Then we have 

(1)
There exists a stability condition $\sigma =(Z_{(\beta, \omega)}, 
\pP)\in \partial U_X$ if and only if $\beta \cdot C \notin \mathbb{Z}$. 
If $\beta \cdot C \in (k-1, k)$, then we have 
$$\pP ((0,1])=\left( \oPPer (X/Y)\otimes \lL ^{\otimes k} \right) \cap 
D(X/S).$$

(2) 
We have 
$$\iPPer (X/Y) \cap D(X/S) =
\left( \oPPer (X/Y) \otimes \lL \right) \cap D(X/S).$$
\end{lem}
\textit{Proof}. 
The proof is same as in~\cite[Lemma 4.3]{Tst},~\cite[Lemma 4.4]{Tst}
and~\cite[Lemma 4.5]{Tst}. $\quad \square$

\vspace{5mm}

Now we can glue the regions $U(W, \Phi)$ 
at the codimension one boundary in both type I and II cases. 
\begin{prop}
The regions $U (W, \Phi)$ satisfy the 
following:
\begin{itemize}
\item
$U(W, \Phi) \cap U(W', \Phi ')\neq \emptyset$ if and only if
$W\cong W'$ and $\Phi ^{-1}\circ \Phi ' \cong \otimes \lL \circ \phi ^{\ast}$
for some $\lL \in \Pic (W)$ and $\phi \in \Aut (W/S)$. In this case, 
we have $U(W, \Phi)=U(W', \Phi ')$. 
\item $U (W, \Phi)\cap U (W', \Phi ')\neq 
\emptyset$ in a codimension one wall of type I
if and only if $W' \dashrightarrow W$ is a flop, 
and $\Phi ^{-1}\circ \Phi ' \cong \otimes 
\lL \circ \phi ^{\ast} \circ \Phi _l \circ \otimes 
\lL ' \circ \phi ^{'\ast}$ 
with $\Phi _l$ a
 standard equivalence and $\lL \in \Pic (W)$, $\lL '\in 
\Pic (W')$, $\phi \in \Aut (W/S)$ and $\phi ' \in \Aut (W'/S)$. 
\item $U (W, \Phi) \cap U (W', \Phi ')\neq 
\emptyset$ in a codimension one wall of type II
if and only if $W\cong W'$ and
$\Phi ^{-1}\circ \Phi ' \cong \otimes 
\lL \circ \phi ^{\ast} \circ \Phi _l \circ \otimes 
\lL ' \circ \phi ^{'\ast}$ 
with $\Phi _l$ an
$EZ$-spherical twist and $\lL \in \Pic (W)$, $\lL '\in 
\Pic (W')$, $\phi \in \Aut (W/S)$ and $\phi ' \in \Aut (W'/S)$. 
\end{itemize}
\end{prop}
\textit{Proof}. The same proof of~\cite[Proposition 4.10]{Tst}
 can be applied. $\quad \square$
 
\subsection*{Descriptions of $\Stabn ^{\circ}(X/S)$}
Before describing $\Stabn ^{\circ}(X/S)$, we give the 
definition of the subset $\FM _{\rm{bir}} ^{\circ}(X)\subset 
\FM _{\rm{bir}}(X)$. 
\begin{defi}
\label{FM}
We define 
$
\FM _{\emph{bir}}^{\circ}(X)\subset  \FM _{\emph{bir}}(X)$ to be the 
subset 
 of pairs $(W,\Phi)\in \FM _{\emph{bir}}(X)$ such that
there exists a sequence of birational maps,
$$W=X^{n}\dashrightarrow X^{n-1}\dashrightarrow \cdots \dashrightarrow 
X^1 \dashrightarrow X^0 =X, $$
and equivalences of birational Fourier-Mukai type over $S$, 
$\Phi ^j \colon D^b (X^j) \to D^b (X^{j-1})$ such that 
$\Phi \cong \Phi ^1 \circ \cdots \circ \Phi ^n$. Each 
$\Phi ^j$ is one of the following:
\begin{itemize}
\item \emph{type I :}  $X^j \dashrightarrow X^{j-1}$ is a flop and 
$\Phi ^j$ is a standard equivalence. 

\item \emph{type I\hspace{-.1em}I :} $X^j = X^{j-1}$ and $\Phi ^j$
is an $EZ$-spherical twist.  

\item \emph{type I\hspace{-.1em}I\hspace{-.1em}I :}
$X^j = X^{j-1}$ and $\Phi ^j$
is a tensoring line bundle $\lL \in \Pic (X^j)$.   

\item \emph{type I\hspace{-.1em}V :}
$X^j =X^{j-1}$ and $\Phi ^j \cong \phi ^{\ast}$ for some 
$\phi \in \Aut (X^j /S)$. 
\end{itemize}
\end{defi}

Now we have the following:
\begin{thm}\label{decom}
We have a disjoint union of locally finite chambers:
$$\mM \cneq 
\bigcup _{(W, \Phi)\in \FM _{\emph{bir}}^{\circ}(X)}
U(W, \Phi) \subset \Stabn ^{\circ}(X/S),$$
in the sense that if two chambers intersect, then they 
coincide. Moreover we have 
$\overline{\mM}=\Stabn ^{\circ}(X/S)$. 
\end{thm}
\textit{Proof}. 
By Lemma 4.5, for $\sigma \in \partial U_X$ there exists a closed point
$x\in X_0$ such that $\oO _x$ is not stable in $\sigma$. Thus $U_X$ is 
one of the connected components of the open subset, 
$$\widetilde{U}_X \cneq \{ \sigma \in \Stabn ^{\circ}(X/S) 
\mid \oO _x \mbox{ is stable for any }x\in X_0 \}.$$
Let us take $\sigma _0 \in U_X$, $\sigma \in \Stabn ^{\circ}(X/S)$,
and a path $\gamma \colon [0,1]\to \Stabn ^{\circ}(X/S)$ such that 
$\gamma (0)=\sigma _0$ and $\gamma (1)=\sigma$. Note that we have the wall 
and chamber structure in the sense of 
Proposition~\ref{wall} by Lemma~\ref{di}.   
One can choose a compact subset $O\subset \Stabn ^{\circ}(X/S)$
 at which $\gamma ((0,1])$ is 
contained in its interior. Then there are finitely many codimension 
one walls in $O$ at which an object $E$ with $[E]=[\oO _x] \in \nN (X/S)$
can become unstable. 

By deforming $\gamma$ a little bit, we may assume 
there exists a finite sequence $0<t_1 <\cdots <t_{n-1}<t_n =1$ 
such that $\gamma (t_{i-1}, t_i)$ is contained in one of the chambers, and 
each $\gamma (t_i)$ is contained in only one wall. Then 
$\gamma ((0, t_1))$ is contained in $U_X$ and $\gamma ((t_1, t_2))$ 
is contained in $U(W, \Phi)$ for some $(W, \Phi)\in 
\FM _{\rm{bir}}^{\circ}(X)$ by 
Lemma~\ref{bou}. Repeating this argument, we can conclude 
$\gamma (t_n)\in \overline{U}(W, \Phi)$ for some $(W, \Phi)\in 
\FM _{\rm{bir}}^{\circ}(X)$. $\quad \square$

\vspace{3mm}

Next we show $\Stabn ^{\circ}(X/S)$ is a regular 
covering space over some open subset of $N^1 (X/S)_{\mathbb{C}}$. 
For the technical reason, we assume the following, 
$$(\star) \quad \mbox{ For a general hyperplane section }
0\in T\subset S, \mbox{ the pull back } X_{T} \cneq f^{-1}(T)
\mbox{ is smooth }.$$
We introduce some notations. 
We define 
$\Lambda _0 \subset \overline{\NE}(X/S)$ to be the 
subset which consists of numerical classes of smooth
rational curves $C\subset X_0$ which 
generate extremal rays of $\overline{\NE}(X/S)$. Also 
define 
$\Lambda \subset \overline{\NE}(X/S)$ to be 
the numerical classes of cycles 
$l=\sum _{i=1}^{k} [C_i]$ with $[C_i]\in \Lambda _0$ and
the dual graph of $C_1 \cdots C_k$ is of Dynkin type. 
Let $E\subset X$ be a $f$-vertical divisor. We define 
$w _E \in \GL (N_1 (X/S))$ to be the reflection
$$w _E (x)\cneq x + (x\cdot E)E_T,$$
for $x\in N_1 (X/S)$. Here $E_T \in N_1 (X/S)$ is the 
fundamental cycle for the scheme theoretic intersection of 
$E$ and $X_T$. We denote by $W_{\rm{ref}}\subset \GL (N_1 (X/S))$ the
subgroup generated by $w _E$ for $f$-vertical divisors $E$. 
We use the following easy lemma.
\begin{lem}\label{eff}
Let us take
 $w \in W_{\emph{ref}}$ and $l\in \Lambda$. Then 
  $w (l)$ or $-w (l)$ is represented by an effective one 
  cycle in $N_1 (X/S)$. 
  \end{lem}
 \textit{Proof}. 
 Let $C_i \subset X_0$, ($i=1, \cdots, N$) 
 be the irreducible components of $X_0$. Note that
 $N_1 (X/S)$ is identified with the vector space 
 $\oplus _{i=1}^N \mathbb{R}C_i$, and we can introduce the 
 bilinear pairing on $N_1 (X/S)$
  by $(C_i, C_j)\cneq (C_i \cdot C_j)_{X_T}$.
 Note that $w _E$ preserves the above bilinear pairing.
 Since $l^2 =-2$ we have $w(l)^2 =-2$. 
 By Zariski's lemma for $X_T$, $w(l)$ or $-w(l)$ is 
 effective. 
 $\quad \square$. 
 
 \vspace{5mm}
 
For $l\in N_1 (X/S)$, we define $H_{l}$ to be the hyperplane:
$$H_{l}\cneq \{ \beta +i\omega \in N^1 (X/S)_{\mathbb{C}}
\mid (\beta +i\omega)\cdot l\in \mathbb{Z} \}.$$
On the other hand, 
let $\Auteq ^{\circ}(X/S)$ be the group of 
autoequivalence of $D(X)$ which is of Fourier-Mukai type over $S$ 
and preserve the connected component $\Stabn^{\circ}(X/S)$. 
We define the group $G$ to be 
$$ 
G  \cneq \ker \left(\Auteq ^{\circ}(X/S)\ni \Phi \longmapsto 
\widetilde{\phi}\in \GL (N^1 (X/S)_{\mathbb{C}}) \right).
$$
Now we have the following:
\begin{thm}\label{norm}
Under the assumption $(\star)$, 
we have the map
$$\zZ _{\rm{n}}\colon \Stabn ^{\circ}(X/S) \lr \bB (X/S)_{\mathbb{C}}
\setminus \bigcup _{(w,l)\in W_{\emph{ref}} \times \Lambda}H_{w(l)}, $$
which is a regular covering map with Galois group equal to $G$. 
\end{thm}
\textit{Proof}. We use the same strategy as in~\cite[Theorem 4.13]{Tst}. 
\begin{step}
We have $\Imm \zZ _{\rm{n}} \subset \bB (X/S)_{\mathbb{C}}\setminus 
\bigcup _{(w,l)\in W_{\emph{ref}} \times \Lambda}H_{w(l)}. $
\end{step}
\textit{Proof}. 
By Proposition~\ref{com} and Lemma~\ref{bound}, we have 
$\zZ _{\rm{n}}(\overline{U}(W, \Phi)) \subset 
\bB (X/S)_{\mathbb{C}}$ for 
$(W, \Phi)\in \FM ^{\circ}_{\rm{bir}}(X)$. 
Therefore by Theorem~\ref{decom}, it suffices to show 
$\zZ _{\rm{n}}(\overline{U}_X)\cap H_{w(l)} =\emptyset$
for $(w,l)\in W_{\rm{ref}} \times 
\Lambda$. Take $\sigma \in \overline{U}_X$, 
$\beta +i\omega \cneq \zZ _{\rm{n}}(\sigma)\in N^1 (X/S)_{\mathbb{C}}$
and assume $(\beta +i\omega)\cdot w(l)\in \mathbb{Z}$ for 
$(w, l)\in W_{\rm{ref}} \times \Lambda$.  
Note that we have $\omega \in \overline{\aA}(X/S) \cap 
\bB (X/S).$
By the base point free theorem~\cite{KMM}, there exists a birational 
contraction over $S$, $X\stackrel{g}{\to}Y \to S$ and an ample 
$\mathbb{R}$-divisor $\omega _Y$ on $Y$ such that $\omega =g^{\ast}
\omega _Y$. 
Note that we have the natural embedding $N_1 (X/Y)\hookrightarrow N_1 (X/S)$ 
and
let $\Lambda ' \cneq \overline{\NE}(X/Y) \cap \Lambda$,  
$\Lambda _0 ' \cneq \overline{\NE}(X/Y) \cap \Lambda _0$. 
Since $w(l)$ or $-w(l)$ is effective by Lemma~\ref{eff}, and $\omega \cdot 
w(l)=0$, 
$w(l)$ is represented by an one cycle contracted by $g$. Also 
$w(l)^2 =-2$ implies $w(l) \in \Lambda '$. Thus we may assume 
$w=\id$ and $l \in \Lambda '$. First assume 
$l\in \Lambda ' _0$ and take a rational curve 
$C\subset X_0$ such that $l=[C]$. 
Then 
$\oO _C (k-1)$ is stable in $U_X$ for $k\in \mathbb{Z}$,
 hence at least semistable 
in $\sigma$. Therefore 
$$\zZ _{(\beta, \omega)}(\oO _{C}(k-1))=
-k +(\beta +i\omega)\cdot C \neq 0.$$ 
Thus $(\beta +i\omega)\cdot l \notin \mathbb{Z}$ 
for $l\in \Lambda ' _0$. 

Next assume $(\beta +i\omega)\cdot l \in \mathbb{Z}$ for some 
$l\in \Lambda '$. 
Then we can find $(W, \Phi)\in \FM^{\circ}(X/Y)$ and an irreducible rational 
curve $C' \subset W$ with $[C']\in \overline{\NE}(W/S)$ an 
extremal ray,  
such that $\widetilde{\phi}(H_{[C']})=H_l$ and $\sigma \in
 \overline{U}(W, \Phi)$. (See the proof of \cite[Theorem 4.13 Step 1]{Tst}.)
 Since 
$\zZ _{\rm{n}}(\overline{U}_W)\cap H_{[C']}=\emptyset$, we have 
$$\beta +i\omega \in \zZ _{\rm{n}}(\overline{U}(W, \Phi))\cap H_l
=\widetilde{\phi}(\zZ _{\rm{n}}(\overline{U}_X) \cap H_{[C']}) =\emptyset$$
by Proposition~\ref{com}, thus a contradiction.  $\quad \square$ 

\begin{step}
The map $\zZ _{\rm{n}}$ is surjective. 
\end{step}
\textit{Proof}. 
Applying $EZ$-spherical twists and flops, 
 it suffices to show the surjectivity on 
$$\left( \overline{\aA}(X/S)_{\mathbb{C}} \cap \bB (X/S)_{\mathbb{C}}
\right) \setminus \bigcup _{(w,l)\in W_{\rm{ref}} \times \Lambda}H_{w(l)}
=\left( \overline{\aA}(X/S)_{\mathbb{C}} \cap \bB (X/S)_{\mathbb{C}}
\right) \setminus \bigcup _{l\in \Lambda}H_l.$$
Take $\beta +i\omega 
\in ( \overline{\aA} (X/S)_{\mathbb{C}}\cap \bB (X/S)_{\mathbb{C}}
)\setminus \cup _{l\in \Lambda}H_l$. Then $\omega$
corresponds to the birational contraction, 
$X \stackrel{g}{\to}Y \to S$, i.e. 
$$ \beta +i\omega 
\in \wW \cneq N^1 (X/S)\oplus ig^{\ast}\aA (Y/S).$$ 
Let $\Lambda _0 '$, $\Lambda '$ be as in Step 1. 
For $l\in \Lambda '$
we define $\overline{H}_l$ to be
$$\overline{H}_l \cneq \{ \beta \in N^1 (X/Y) \mid \beta \cdot l \in 
\mathbb{Z} \} \subset N^1 (X/Y).$$
Note that $\overline{H}_l$ is a real codimension one hypersurface
in $N^1 (X/Y)$. 
The composition $\wW \to N^1 (X/S) \to N^1 (X/Y)$ induces the 
following topological fiber space structure:
$$\Pi \colon \wW \setminus \bigcup _{l\in \Lambda}H_l \lr N^1 (X/Y) \setminus 
\bigcup _{l\in \Lambda '}\overline{H}_l.$$
Let $C(X/Y)$ be one of the connected component of the right hand side,
$$C(X/Y)\cneq \{ \beta \in N^1 (X/Y) \mid -1 < \sum _{l\in \Lambda _0 '} \beta \cdot l <0, -1< \beta \cdot l <0 \mbox{ for all }
l\in \Lambda _0 ' \}.$$
Then by the argument of~\cite[Theorem 4.13, Step 2]{Tst}, 
we can find $(W, \Phi)\in \FM^{\circ}(X/Y)$ such that $\widetilde{\phi}^{-1}
\Pi (\beta +i\omega)\in N^1 (W/Y)$ is contained in $C(W/Y)$. Thus we may assume 
$\Pi (\beta +i\omega) \in C(X/Y)$ by Proposition~\ref{com}. 
In this region, the pair 
$$\sigma \cneq (Z_{(\beta, \omega)}, \oPPer (X/Y) \cap D(X/S)), $$
gives a point of $\Stabn ^{\circ}(X/S)$.
In fact let $D(X/Y)\subset D(X/S)$ be the subcategory whose objects 
are supported on $\Ex (g)\cap X_0$. Then
the pair $(Z_{(\beta, \omega)}, \oPPer (X/Y) \cap D(X/Y))$ 
gives a point of $\Stab (D(X/Y))$ by~\cite[Theorem 4.13, Step 2]{Tst}.
(Note that we assumed $g \colon X\to Y$ to be small in~\cite{Tst}. 
But we can easily generalize the argument of~\cite{Tst} for the
case of birational contraction $g \colon X\to Y$, when dimensions of all the 
fibers are less than or equal to zero.)
Thus $Z_{(\beta, \omega)}(E)$ is contained in the upper half plane 
if $E\in \oPPer (X/Y)\cap D(X/Y)$. 
Also if $E\in \oPPer (X/Y) \cap D(X/S)$ and not contained in 
$D(X/Y)$, then $\Imm Z_{(\beta, \omega)}(E)>0$. The Harder-Narasimhan 
property is satisfied by the same argument of~\cite[Lemma 4.5]{Tst}. 
$\quad \square$

\begin{step}
$\zZ _{\rm{n}}$ is a regular covering map with Galois group equal to $G$. 
\end{step}
\textit{Proof}. 
Note that $G$ acts on $\Stabn ^{\circ}(X/S)$ as deck transformations. 
Let us take $\sigma, \sigma ' \in \Stab ^{\circ}(X/S)$ such that 
$\zZ _{\rm{n}}(\sigma)=\zZ _{\rm{n}}(\sigma ')$. We will find 
$g\in G$ such that $g(\sigma)=\sigma '$. 
By Theorem~\ref{decom}, we may assume $\sigma \in U_X$ and
 $\sigma ' \in U(X, \Phi)$ for some $\Phi \in \Auteq ^{\circ}(X/S)$. 
Then $\widetilde{\phi}\in \GL (N^1 (X/S)_{\mathbb{C}})$
 is written as the composition 
$\widetilde{\phi}=w \circ \phi _{\ast}\circ \otimes \lL$.  
Here $w$ is a composition of reflections with respect to 
$f$-vertical divisors, $\phi$ is a birational map given by $\Phi$ and
$\lL \in \Pic (X)$. Since $\zZ _{\rm{n}}(\sigma)=\zZ _{\rm{n}}(\sigma ')$, 
we have 
$$\widetilde{\phi}(\aA (X/S)_{\mathbb{C}}) \cap \aA (X/S)_{\mathbb{C}}
\neq \emptyset. $$
Since $\otimes \lL$ preserves $\aA (X/S)_{\mathbb{C}}$ and 
$\phi _{\ast}$ preserves $\overline{\mM}(X/S)_{\mathbb{C}}$, 
we have 
$$w (\overline{\mM} (X/S)^{\circ}_{\mathbb{C}}\cap \bB (X/S)_{\mathbb{C}}) 
\cap (\overline{\mM} (X/S)^{\circ}_{\mathbb{C}}\cap \bB (X/S)_{\mathbb{C}}) 
\neq \emptyset.$$
Here $\overline{\mM}(X/S)_{\mathbb{C}}^{\circ}$ is the set of 
inner points of $\overline{\mM}(X/S)_{\mathbb{C}}$. 
Since $w$ is a composition of reflections, it follows that 
$w =\id$. Therefore $\phi _{\ast}\aA (X/S)_{\mathbb{C}}\cap 
\aA (X/S)_{\mathbb{C}}\neq \emptyset$, and this implies 
$\phi$ gives a $S$-isomorphism by Lemma~\ref{mov}. 
Let $g\cneq \Phi \circ \phi _{\ast}^{-1}\circ \otimes \lL ^{-1}\in 
\Aut ^{\circ}(X/S)$. 
Then $g \in G$, $g(\sigma) \in U(X, \Phi)$ and $\zZ _{\rm{n}}\circ g(\sigma)
=\zZ _{\rm{n}}(\sigma)$. Therefore $\sigma '=g(\sigma)$. $\quad \square$

\subsection*{Non-normalized stability conditions}
We give the description of $\Stab ^{\circ}(X/S)$. Note that 
we have the natural map
$$\alpha \colon \mathbb{C}\times \Stabn ^{\circ}(X/S) 
\lr \Stab ^{\circ}(X/S),$$
given by the rescaling action of $\mathbb{C}$ on $\Stab (X/S)$. We 
have the following proposition:
\begin{prop}\label{non-nor}
$\alpha$ is an isomorphism and we have the following commutative diagram:
$$\xymatrix{
\alpha \colon \mathbb{C}\times \Stabn ^{\circ}(X/S) \ar[r]^{\cong}\ar[d]_{1 \times \zZ _{\rm{n}}} & 
\Stab ^{\circ}(X/S) \ar[d]^{\zZ} \\
e \colon \mathbb{C}\oplus N^1 (X/S)_{\mathbb{C}} \ar[r] & \mathbb{C}
\oplus N^1 (X/S)_{\mathbb{C}}. }$$
Here $e$ takes $(\lambda, L)$ to 
$(\exp (-i\pi \lambda), \exp (-i\pi \lambda)L)$. 
\end{prop}
\textit{Proof}. 
The same proof of~\cite[Theorem 5.5]{Tst} works. $\quad \square$

\vspace{5mm}
By Theorem~\ref{decom} and Proposition~\ref{non-nor}, 
$\Stab ^{\circ}(X/S)$ is a regular covering space over
$$ e( \mathbb{C} \oplus 
( \bB (X/S)_{\mathbb{C}}\setminus \cup _{(w,l)\in W_{\rm{ref}}
\times \Lambda}H_{w(l)} ))
\subset \mathbb{C}\oplus N^1 (X/S)_{\mathbb{C}}.$$
We give the explicit description of the above set.   
Let $V\subset N^1 (X/S)$ be the subspace generated by 
$f$-vertical divisors. 
If we choose $H\in N^1 (X/S)$ to be $\deg (H|_{X_{\bar{\eta}}})=1$, 
we have the decomposition, 
$$\mathbb{C}\oplus 
N^1 (X/S)_{\mathbb{C}} \ni (\lambda, D) \longmapsto (\lambda, 
\deg(D|_{X_{\bar{\eta}}}), D-\deg 
(D|_{X_{\bar{\eta}}})H) \in \mathbb{C}^{2} \oplus V_{\mathbb{C}}.$$
Under the above decomposition, the subset $e(\mathbb{C}\oplus 
\bB (X/S)_{\mathbb{C}}) \subset \mathbb{C}\oplus N^1 (X/S)_{\mathbb{C}}$
corresponds to $\GL ^{+}(2, \mathbb{R}) \times V_{\mathbb{C}}$. 
Here $\GL ^{+}(2, \mathbb{R})$ is a subset of $\GL (2, \mathbb{R})$ 
which preserves the orientation, and embedded into $\mathbb{C}^2$ 
via 
$$\GL ^{+}(2, \mathbb{R}) \ni \left( \begin{array}{cc}
a & b \\ c & d \end{array} \right)  \longmapsto 
(a+ci, b+di)\in \mathbb{C}^2.$$
For $k\in \mathbb{Z}$ and $l\in N_1 (X/S)$, define 
$\widetilde{H}_{k,l}$ to be the codimension two hyper plane, 
$$\widetilde{H}_{k,l}\cneq 
\left\{ \left( \left( \begin{array}{cc} a & b \\ c & d \end{array} \right), 
v_0 +v_1 i \right) 
\in \GL^{+}(2, \mathbb{R}) \times V_{\mathbb{C}}: \begin{array}{ll}
v_0 \cdot l + bH\cdot l -ka =0 \\
v_1 \cdot l +dH\cdot l -kc =0 
\end{array} \right\}.$$
Then the image of 
$\mathbb{C}\oplus H_l \subset \mathbb{C}\oplus N^1 (X/S)_{\mathbb{C}}$
under $e$ is $\cup _{k\in \mathbb{Z}}\widetilde{H}_{k,l}$. 
Thus we have the following: 

\begin{thm}
Under the assumption $(\star)$, 
there exists 
a map 
$$\zZ \colon \Stab ^{\circ}(X/S) \lr (\GL ^{+}(2, \mathbb{R}) \times 
V_{\mathbb{C}})\setminus \bigcup _{(k,w,l)\in 
\mathbb{Z}\times W_{\emph{ref}} \times \Lambda}\widetilde{H}_{k, w(l)},$$
which is a regular covering map.
\end{thm}
\section{Localization theorem for stability conditions}
In this section we 
establish the property of stability conditions,
 which is similar to the 
localization theorem for Grothendieck groups. 
Let $f\colon X\to S$ be a smooth 
projective morphism with $S=\Spec \mathbb{C}[[t]]$. 
Here 
we don't have to assume $f$ is a Calabi-Yau fibration. 
Let $X_0$ be the closed fiber of $f$ and $i\colon X_0 \hookrightarrow 
X$ be the inclusion. Also let $\Stab (X_0)$, $\Stab (X/S)$ be
the spaces of numerical stability conditions as in Section 2.
 \begin{prop}\label{mo}
There exists a map $\theta \colon \Stab (X/S) \to 
\Stab (X_0)$ which fits into the diagram,  
$$\xymatrix{
\Stab (X/S) \ar[d]_{\zZ} \ar[r] ^{\theta} &  \Stab (X_0) \ar[d]^{\zZ _0} \\
\nN (X)_{\mathbb{C}} \ar[r]^{i^{\ast}} & \nN (X_0)_{\mathbb{C}}.
}$$ 
\end{prop}
\textit{Proof}. The construction of $\theta$ is due 
to~\cite[Corollary 2.2.2]{Pol}. According to \textit{loc.cit.}, 
for $\sigma =(Z, \pP)\in \Stab (X/S)$, $\theta (\sigma)$ can be constructed
to be the pair $(Z_0, \pP _0)$ with $Z_0 =i^{\ast}Z$ and 
$$\pP _0 (\phi) =\{ E\in D(X_0) \mid i_{\ast}E \in \pP (\phi) \}.$$
In \textit{loc.cit.}, the assumption $\oO _{X_0}\dotimes \pP (t) 
\subset \pP (t, \infty)$ was needed. But this assumption is 
satisfied in our case, as proved in~\cite[Theorem 2.3.5]{Pol}. 
$\quad \square$.

\vspace{5mm}

Note that $i^{\ast}\colon \nN (X)_{\mathbb{C}} \to \nN (X_0)_{\mathbb{C}}$ is
injective since $i_{\ast}\colon K(X_0) \to K(X/S)$ is an isomorphism. 
The map $\theta \colon \Stab (X/S) \to \Stab (X_0)$ induces the 
map 
$$\widetilde{\theta} \colon \Stab (X/S) \lr \widetilde{\Stab}(X_0) \cneq 
\Stab (X_0) \times _{\nN (X_0)_{\mathbb{C}}}\nN (X)_{\mathbb{C}}.$$
The purpose of this section is to study the above map. 
We prepare some lemmas.

\begin{lem}\label{iso}
Let us take $A, B\in D(X_0)$ which satisfy $\Ext _{X_0}^{i}(A, B)=0$ 
for $i=-1,-2$. Then $i_{\ast}\colon \Hom _{X_0}(A,B)\to \Hom _{X}(i_{\ast}A, 
i_{\ast}B)$ is an isomorphism. 
\end{lem}
\textit{Proof}.
Note that we have $\Hom _{X}(i_{\ast}A, i_{\ast}B)\cong 
\Hom _{X_0}(\dL i^{\ast}i_{\ast}A, B)$ by adjunction. 
We have the distinguished triangle in $D(X_0)$, 
$$A[1] \otimes \oO _{X_0}(X_0)=A[1] \lr \dL i^{\ast}i_{\ast}A \lr A \lr A[2],$$
by~\cite[Lemma 3.3]{B-O2}. Applying $\Hom _{X_0}(\ast, B)$, we 
obtain the long exact sequence:
$$\Ext _{X_0}^{-2}(A, B) \lr \Hom _{X_0}(A, B) \lr 
\Hom _{X_0}(\dL i^{\ast}i_{\ast}A, B) \lr \Ext _{X_0}^{-1}(A,B).$$
Since $\Ext _{X_0}^{-2}(A,B)=\Ext _{X_0}^{-1}(A,B)=0$, the 
morphism $i_{\ast}\colon \Hom _{X_0}(A,B)\to \Hom _{X}(i_{\ast}A, 
i_{\ast}B)$ is an isomorphism. $\quad \square$

\begin{lem}\label{ff}
Take $\sigma =(Z, \pP)\in \Stab (X/S)$ and let
$\theta (\sigma)=(Z_0, \pP _0)\in \Stab (X_0)$. Let
$\aA _0 \cneq \pP _0 ((0,1])$ and $\aA \cneq \pP ((0,1])$ 
be the corresponding Abelian categories. 
Then we have 

(i) $i_{\ast}\colon \aA _0 \to \aA$ is fully faithful and has 
a left adjoint $i^{\ast}\colon \aA \to \aA _0$. 

(ii) If $E\in \aA$ satisfies $\Hom (E, E)=\mathbb{C}$, 
 then $E\cong i_{\ast}i^{\ast}E$. 
 
(iii) We have $i_{\ast}\pP _{0,s}(\phi)=\pP _{s}(\phi)$ for 
$\phi \in \mathbb{R}$. 

(iv) For $E\in D(X_0)$ 
and $U\in \nN (X)_{\mathbb{C}}$, we have 
$\phi ^{\pm}_{\sigma}(i_{\ast}E)=\phi ^{\pm}_{\theta (\sigma)}(E)$ and
$\lVert U \rVert _{\sigma}
=\lVert i^{\ast}U \rVert _{\theta (\sigma)}$.

(v) For another $\tau =(W, \qQ)\in \Stab (X/S)$, we have 
$d(\sigma, \tau)=d(\theta (\sigma), \theta (\tau))$. 
\end{lem}
\textit{Proof}. 
(i) By lemma~\ref{iso}, $i_{\ast}\colon \aA _0 \to \aA$ is fully faithful. 
We define $i^{\ast}\colon \aA \to \aA _0$ to be 
$$i^{\ast}(F) \cneq H_{\aA _0}^{0}(\dL i^{\ast}F)\in \aA _0.$$
Here $H_{\aA _0}^0(\ast)$ is the zero-th cohomology functor
for the t-structure on $D(X_0)$ with heart $\aA _0$. 
We check $i^{\ast} \colon \aA \to \aA _0$ gives a left adjoint of 
$i_{\ast}\colon \aA _0 \to \aA$. Take $A\in \aA$ and $B\in \aA _0$. 
Then we have 
\begin{align*}
\Hom _{X}(A, i_{\ast}B) & \cong \Hom _{X_0}(\dL i^{\ast}A, B) \\
& \cong \Hom _{X_0}(\tau _{\ge 0}\dL i^{\ast}A, B).
\end{align*}
Here $\tau _{\ge 0}$ is a truncation functor for
the t-structure with heart $\aA _0$. Thus it suffices to 
check $H_{\aA _0}^{j}(\dL i^{\ast}A)=0$ for $j\ge 1$. 
Assume the contrary. Then there exists $j\ge 1$ and 
a non-zero map $\dL i^{\ast}A \to C[-j]$
for some $C\in \aA _0$. Taking the adjoint, we obtain the non-zero map
$A\to i_{\ast}C[-j]$. Since $i_{\ast}C\in \aA$, this is a 
contradiction. 

(ii) Assume $E\in \aA$ satisfies $\Hom (E, E)=\mathbb{C}$. We have 
\begin{align*}
i_{\ast}i^{\ast}E &= i_{\ast}H_{\aA _0}^0 (\dL i^{\ast}E) \\
&\cong H_{\aA}^0 (i_{\ast}\dL i^{\ast}E) \\
& \cong H_{\aA}^0 (E\dotimes _{\oO _X}\oO _{X_0}).
\end{align*}
The second isomorphism follows from the fact that 
 $i_{\ast}$ takes $\aA _0$ to 
$\aA$ and the third isomorphism follows from the projection formula. 
We have the exact sequence:
$$0 \lr \oO _X \stackrel{\times t}{\lr} \oO _X \lr \oO _{X_0} \lr 0.$$
Applying $E\dotimes _{\oO _X}$, we obtain the distinguished triangle
$$E \stackrel{\times t}{\lr} E \lr E\dotimes _{\oO _X}\oO _{X_0}
\lr E[1].$$
Since $\Hom _{X}(E,E)=\mathbb{C}$, the map
$E\stackrel{\times t}{\to}E$ is an isomorphism or zero-map. 
Let us assume $E\stackrel{\times t}{\to}E$
 is an isomorphism. Then $E\stackrel{\times t^n}{\to}E$ also 
gives an isomorphism for $n>0$. But $H^i(E) \stackrel{\times t^n}{\to}
H^i(E)$ must be zero map for some $n>0$ since $\Supp H^i (E)\subset X_0$. 
This is a
contradiction and it follows that 
$E\stackrel{\times t}{\to}E$ must be zero map. 
Thus we have the decomposition
$$E \dotimes _{\oO _X}\oO _{X_0} \cong E \oplus E[1], $$
and we have $H_{\aA}^0 (E\dotimes _{\oO _X}\oO _{X_0}) \cong E$.

(iii)
Take $E\in \pP _{s}(\phi)$ such that $0<\phi \le 1$. 
We denote by $\phi _0 (\ast)$, $\phi (\ast)$ the phases of 
objects in $\aA _0$, $\aA$ for stability functions 
$Z_0$, $Z$ respectively. 
By (i), $i_{\ast}\colon \aA _0 \to \aA$ has a left adjoint 
$i^{\ast}\colon \aA \to \aA _0$, and
since $\Hom (E, E)=\mathbb{C}$, we have $E\cong i_{\ast}i^{\ast}E$.
 We show $i^{\ast}E\in \aA _0$ is stable in $\sigma _0$. 
 Assume the contrary. 
Then there exists a monomorphism $F\hookrightarrow i^{\ast}E$ 
in $\aA _0$ such that $\phi _0 (F)\ge \phi _0 (i^{\ast}E)$. 
Since $i_{\ast}$ takes $\aA _0$ to $\aA$, we have the monomorphism 
$i_{\ast}F \hookrightarrow i_{\ast}i^{\ast}E \cong E$
such that $\phi (i_{\ast}F)\ge \phi (E)$. This contradicts that 
$E$ is stable. Therefore $i^{\ast}E\in \aA _0$ is stable in $\sigma _0$, 
and $E\cong i_{\ast}i^{\ast}E$ implies 
$i_{\ast}\pP _{0,s}(\phi)\supset \pP _s (\phi)$. 
Conversely take a stable object $E\in \aA _0$ and assume $i_{\ast}E \in \aA$ 
is not stable. Then there exists a stable object $F\in \aA$ such that 
there exists a monomorphism $F\hookrightarrow i_{\ast}E$
and $\phi (F)\ge \phi (i_{\ast}E)$. Since $i_{\ast}i^{\ast}F\cong F$ and 
$i_{\ast}\colon \aA _0 \to \aA$ is fully faithful, we obtain the 
monomorphism $i^{\ast}F \hookrightarrow E$ such that $\phi _0 (i^{\ast}F)
\ge \phi _0 (E)$. But this also contradicts that $E$ is stable, 
thus we have $i_{\ast}\pP _{s,0}(\phi)\subset \pP _s (\phi)$. 
Consequently we have $i_{\ast}\pP _{0,s}(\phi)=\pP _{s}(\phi)$, which 
is (iii). 

(iv)
This follows from (iii) and the definitions of 
$\phi _{\sigma}^{\pm}$, $\phi _{\theta(\sigma)}^{\pm}$ and 
$\lVert \ast \rVert _{\sigma}$,
$\lVert \ast \rVert _{\theta (\sigma)}$. 

(v) By the definition of $d(\ast, \ast)$ and (iv), 
we have $d(\theta (\sigma), \theta (\tau))\le d(\sigma, \tau)$. 
We show the converse inequality. For instance denote $\varepsilon =
d(\theta (\sigma), \theta (\tau))$, $\theta (\tau)=(i^{\ast}W, \qQ _0)$. 
Let us take $E\in D(X/S)$ and let $A_j \in \pP _s (\phi _j)$ be 
stable factors in $\sigma$ such that 
$\phi _1 = \phi _{\sigma}^{+}(E)$ and 
$\phi _n = \phi _{\sigma}^{-}(E)$. 
Then $i^{\ast}A_j \in \pP _{0,s}(\phi _j)$, thus 
$i^{\ast}A_j \in \qQ _0 ([\phi _j -\varepsilon, \phi _j +\varepsilon])$. 
Therefore 
$$A_j \cong i_{\ast}i^{\ast}A_j \in \qQ ([\phi _j -\varepsilon, \phi _j +
\varepsilon ]).$$
As a consequence, we have $E\in \qQ ([\phi _n -\varepsilon, \phi _1 +\varepsilon])$. This implies $d(\sigma, \tau)\le \varepsilon$. $\quad \square$

\vspace{5mm}

\begin{cor}\label{con}
The map $\theta \colon \Stab (X/S) \lr \Stab (X_0)$ is 
continuous and injective.  \end{cor}
\textit{Proof}. 
Lemma~\ref{ff} (iv), (v) imply 
 $\theta (B_{\varepsilon}(\sigma)) 
\subset B_{\varepsilon}(\theta (\sigma))$ for
$\sigma \in \Stab (X/S)$ and $0< \varepsilon \ll 1$. 
This implies $\theta$ is continuous. 
We check $\theta$ is injective. In fact assume
$\sigma _i =(Z_i, \pP _i) \in \Stab (X/S)$ for $i=1,2$ 
satisfy $\theta (\sigma _1)=\theta (\sigma _2)$. 
Then Lemma~\ref{ff} (iii) implies $\pP _{1,s}(\phi)=\pP _{2,s}(\phi)$ 
for $\phi \in \mathbb{R}$, thus $\pP _1 (\phi)=\pP _2 (\phi)$. 
Since $i^{\ast}\colon \nN (X)_{\mathbb{C}}\to \nN (X_0)_{\mathbb{C}}$
is injective, one has $Z_1 =Z_2$. Thus $\sigma _1 =\sigma _2$ follows. 
$\quad \square$

\vspace{5mm}

Let us take a connected component $\Stab ^{\circ}(X/S)\subset 
\Stab (X/S)$. Since $\theta$ is continuous, one can find 
a connected component $\widetilde{\Stab} ^{\circ}(X_0)\subset 
\widetilde{\Stab} (X_0)$ such that $\widetilde{\theta}$
 takes $\Stab ^{\circ}(X/S)$
into $\widetilde{\Stab} ^{\circ}(X_0)$. 
The following is the main theorem of this section.

\begin{thm}\label{loc}
Assume that for $\sigma \in \Stab ^{\circ}(X/S)$ the linear subspace 
$\{ U\in \nN (X)_{\mathbb{C}} \mid \lVert U \rVert _{\sigma} <\infty \}
\subset \nN (X)_{\mathbb{C}}$ is defined over $\mathbb{Q}$. 
Then the map
$$\theta ^{\circ}\cneq \widetilde{\theta}|_{\Stab^{\circ}(X/S)}
\colon \Stab ^{\circ}(X/S) \lr \widetilde{\Stab} ^{\circ}(X_0),$$
is a homeomorphism. 
\end{thm}
\textit{Proof}. 
By Corollary~\ref{con}, it remains
 to check the surjectivity of $\theta ^{\circ}$. 
Since $\theta$ satisfies $\zZ _0 \circ \theta =i^{\ast}\circ \zZ$, 
the map $\theta ^{\circ}$ is an open map. Thus it suffices to 
show $\Imm (\theta ^{\circ})\subset \widetilde{\Stab}^{\circ}(X_0)$
is closed. This is equivalent to $\theta (\Stab ^{\circ}(X/S))\subset 
\Stab (X_0)$ is closed. 
Take $\sigma _n =(Z_n, \pP _n)\in \Stab ^{\circ}(X/S)$ 
such that $\theta (\sigma _n)$ converges 
to $\tau _0=(i^{\ast}W, \qQ_0)\in \Stab (X_0)$
for some $W\in \nN (X)_{\mathbb{C}}$.  
We show $\tau _0 \in \theta (\Stab ^{\circ}(X/S))$. 
By the assumption, 
we may assume $Z_n \in \nN (X)_{\mathbb{C}}$ 
are rational for all $n$. 
Fix $0 <\varepsilon <1/16$ which satisfies $2(1+\tan \pi \varepsilon)\tan \pi \varepsilon <\sin \pi /8$. 
 Then there exists $N>0$ such that 
if $n>N$ then $\theta(\sigma _n) \in B_{\varepsilon}(\tau _0)$.
Below we fix such $n>N$. 
 By~\cite[Lemma 6.2]{Brs1}, one can choose a constant $k>0$ which 
depends only $\varepsilon$ such that 
$\lVert U \rVert _{\tau '} <k\lVert U \rVert _{\tau _0}$ for 
every $\tau ' \in B_{\varepsilon}(\tau _0)$ and $U\in \nN (X_0)_{\mathbb{C}}$. 
In fact according to the proof of~\cite[Lemma 6.2]{Brs1}, one
can take $k$ to be $k=(1+\tan \pi \varepsilon)/\cos \pi \varepsilon$. 
Therefore we have 
$$\lVert i^{\ast}W -i^{\ast}Z_n  \rVert _{\theta(\sigma _n)} 
<(1+\tan \pi \varepsilon)\tan \pi \varepsilon .$$
From Lemma~\ref{ff} (iv), we have   
$$\lVert W-Z_n \rVert _{\sigma _n} = \lVert i^{\ast}W -i^{\ast}Z_n 
 \rVert _{\theta(\sigma _n)}
 < \sin \pi \varepsilon ' < \sin \pi /8.$$
 Here we have taken $0<\varepsilon ' <1/8$ to be 
 $$\sin \pi \varepsilon ' = 2(1+\tan \pi \varepsilon)\tan \pi \varepsilon.$$
 By Theorem~\ref{def}, one can construct 
$\tau =(W, \qQ) \in \Stab ^{\circ}(X/S)$ such that  
$d(\tau, \sigma _{n})<\varepsilon '$
uniquely. It is enough to check $\theta (\tau)=\tau _0$. 
Let $\varepsilon '' \cneq \max \{\varepsilon ', 2\varepsilon \}<1/8.$
For $m>N$, we have 
\begin{align*}
d(\sigma _n, \sigma _m) & = d(\theta (\sigma_n), \theta (\sigma _m)) \\
& \le d(\theta (\sigma _n),\tau _0)+d(\tau _0, \theta (\sigma _m)) \\
& <2\varepsilon \le \varepsilon ''
\end{align*}
and 
\begin{align*}
\lVert Z_n -Z_m \rVert _{\sigma _n} &= \lVert i^{\ast}Z_n -i^{\ast}Z_m 
\rVert _{\theta (\sigma _n)} \\
&< k\lVert i^{\ast}Z_n -i^{\ast}Z_m \rVert _{\tau _0} \\
& \le k\lVert i^{\ast}Z_n -W_0 \rVert _{\tau _0}+k\lVert W_0 -i^{\ast}Z_m
\rVert _{\tau _0} \\
& <2(1+\tan \pi \varepsilon)\tan \pi \varepsilon =
 \sin \pi \varepsilon ' \le \sin \pi \varepsilon ''.
\end{align*}
Therefore for $m>N$, we have $\sigma _m \in B_{\varepsilon ''}(\sigma _n)$
and $\tau \in B_{\varepsilon '}(\sigma _n)\subset 
B_{\varepsilon ''}(\sigma _n)$.
Because $\varepsilon ''<1/8$, 
$\zZ |_{B_{\varepsilon ''}(\sigma _n)}$ is homeomorphism onto 
its image by Theorem~\ref{def}. Since $Z_{m}$ converge to $W$, 
 $\sigma _m$ must converge to $\tau$, hence $\theta (\tau)=\tau _0$
 by the continuity of $\theta$. 
$\quad \square$.

\section{Stability conditions for smooth $K3$ (Abelian) fibrations}
Let $f\colon X\to S=\Spec \mathbb{C}[[t]]$ be a three dimensional smooth
Calabi-Yau fibration. Then $X_0$ is a $K3$ surface 
or an Abelian surface. 
Now we can describe the spaces of stability conditions 
on $D(X/S)$, using the result of the last section
and the description of $\Stab (X_0)$ given in~\cite{Brs2}. 
We give the description when $X_0$ is a $K$3 surface. 
The other case is similarly discussed. 
\subsection*{Stability conditions on $K3$ surfaces}
Let us recall the construction of stability conditions on 
$D(X_0)$ via tilting. 

A pair $(\tT, \fF)$ of full subcategories
 of an Abelian category $\aA$ is called 
a torsion pair if $\Hom (\tT, \fF)=0$ and
any object $A\in \aA$, fits into an exact sequence, 
$$0 \longrightarrow T \longrightarrow A \longrightarrow F
 \longrightarrow 0,$$
 with $T\in \tT$ and $F\in \fF$. 
 Given a torsion pair of $\aA$, we can produce
another Abelian category in $\aA ^{\sharp}\subset D^b(\aA)$
to be 
$$\aA ^{\sharp}\cneq \{E\in D^b(\aA) \mid H^i (E)=0 \mbox{ for }i
 \neq \{0,-1\}, H^0 (E)\in \tT, H^{-1}(E)\in \fF\}. $$
 In fact $\aA ^{\sharp}$
 is a heart of some bounded t-structure on $D^b(\aA)$. 
 We say $\aA ^{\sharp}$ is a tilting
  with respect to the torsion pair $(\tT, \fF)$.

 For a torsion free sheaf $E \in \Coh (X_0)$ and 
$\omega \in \aA (X_0)$, let $\mu _{\omega}(E)$ be the slope, 
$$\mu _{\omega} (E)\cneq \frac{c_1 (E)\cdot \omega}{r(E)}.$$
Here $r(E)$ is a rank of $E$. 
One has the Harder-Narasimhan filtration 
$$0= E_0 \subset E_1 \subset \cdots \subset E_{n-1}\subset E_n =E,$$
with $F_i =E_i /E_{i+1}$ is $\mu _{\omega}$-semistable and 
$\mu (F_{i})> \mu (F_{i+1})$. 
For $\beta +i\omega \in \aA (X_0)_{\mathbb{C}}$, 
we define $\tT \subset \Coh (X_0)$ to be
the subcategory consists of sheaves whose torsion free parts 
have $\mu _{\omega}$-semistable Harder-Narasimhan factors 
of slope $\mu _{\omega} > \beta \cdot \omega$. Also 
define $\fF \subset \Coh (X_0)$ to be the subcategory consists of 
torsion free sheaves whose $\mu _{\omega}$-semistable factors have slope 
$\mu _{\omega}\le \beta \cdot \omega$. 
Then the pair $(\tT, \fF)$ is a torsion pair,
 and its 
tilting gives the Abelian category 
$\aA _{(\beta, \omega)}\subset D(X_0)$. 
Then let $Z_{(\beta, \omega)}\colon \nN (X_0)\to \mathbb{C}$ be
$$Z_{(\beta, \omega)}(E) \cneq -\int e^{-(\beta +i\omega)}\ch (E)
 \sqrt{\td _{X_0}}.$$
\begin{lem}\emph{\bf{\cite[Proposition 9.2]{Brs2}}}
For $\beta +i\omega \in \aA (X_0)_{\mathbb{C}}$, 
the pair $\sigma _{(\beta, \omega)}=
(Z_{(\beta, \omega)}, \aA _{(\beta, \omega)})$ 
gives a stability condition on $D(X_0)$ if and only if 
for all spherical sheaves $E$ on $X_0$ one has 
$Z_{(\beta, \omega)}(E) \notin \mathbb{R}_{\le 0}$. 
\end{lem}
The stability conditions $\sigma _{(\beta, \omega)}$ are contained in 
one of the connected component, denoted by $\Stab ^{\circ}(X_0)$. 

Next recall that we have the isomorphism, 
$$\ch(\ast)\sqrt{\td _{X_{0}}} \colon \nN (X_0) \stackrel{\cong}{\lr}
\mathbb{Z} \oplus \NS (X_0) \oplus \mathbb{Z}.$$
Under the above isomorphism, the pairing $-\chi (\ast, \ast)$ 
on the left hand side corresponds to the Mukai bilinear pairing 
$$(r, l, s)\cdot (r', l', s')=l\cdot l' -rs' -r's, $$
on the right hand side. Under the above pairing, we define 
$\pP ^{\pm}(X_0)$ as in~\cite{Brs2}, 
$$\pP ^{\pm}(X_0) \cneq 
\{ v_0 +iv_1 \in \nN (X_0)_{\mathbb{C}} \mid 
v_0, v_1 \mbox{ span a positive definite two plane in }
\nN (X_0)_{\mathbb{R}} \}.$$
Note that $\pP ^{\pm}(X_0)$ consists of two connected components. 
We define $\pP ^{+}(X_0)$ to be one of the connected component 
of $\pP ^{\pm}(X_0)$, which contains $(1, i \omega, -\frac{1}{2}\omega ^2)$.
Let $\Delta (X)$ be
$$\Delta (X) \cneq \{ \delta \in \nN (X_0) \mid \delta ^2 =-2 \},$$
and $\pP _0 ^{+}(X_0)$ be
$$\pP _0 ^{+}(X_0) \cneq \pP ^{+}(X_0) \setminus \bigcup _{\delta 
\in \Delta(X)} \delta ^{\perp}.$$
Here $\delta ^{\perp}\cneq \{ u\in \nN (X_0)_{\mathbb{C}}\mid 
(u, \delta)=0\}.$
\begin{thm}\emph{\bf{\cite[Theorem 1.1]{Brs2}}} \label{K3}
Sending stability conditions to their central charges gives the map 
$$\zZ _0 \colon \Stab ^{\circ}(X_0) \lr \pP _0 ^{+}(X_0), $$
which is a regular covering map. 
\end{thm}
\subsection*{The description of $\Stab (X/S)$.}
Now we can apply theorem~\ref{loc} to give the description of 
$\Stab (X/S)$ when $f\colon X\to S=\Spec \mathbb{C}[[t]]$
 is a three dimensional smooth $K3$ fibration. 
 Let us take $\beta '+i\omega '\in \aA (X/S)_{\mathbb{C}}$ 
 and the restriction to $X_0$, $\beta  +i\omega  \cneq 
 i^{\ast}(\beta '+i\omega')\in \aA (X_0)_{\mathbb{C}}$. 
 Then we can construct the torsion pair
  $(\tT, \fF)$
 on $\Coh (X_0)$ with respect to $\beta +i\omega \in 
 \aA (X_0)_{\mathbb{C}}$, and let 
 $\aA _{(\beta, \omega)}$ be the tilting as before. 
 Now we define the categories 
 $\tT '$ and $\fF '$ to be the minimum extension 
 closed subcategory of $\Coh (X/S)$, which contain $i_{\ast}\tT$ and
 $i_{\ast}\fF$ respectively.
 \begin{lem}
 The pair $(\tT ', \fF ')$ is a torsion pair on $\Coh (X/S)$. 
 \end{lem}
 \textit{Proof}. 
 First $\Hom _{X}(\tT ', \fF')=0$ follows from Lemma~\ref{iso}. 
 Let us take $E\in \Coh (X/S)$. 
 By taking Harder-Narasimhan filtration and  
 Jordan-H\"{o}rder filtration in $\omega '$-Giesker stability, 
 we have the filtration 
 $$0=E_0 \subset E_1 \subset \cdots \subset E_{n-1} \subset E_n =E, $$
 such that each quotient $F_j =E_j /E_{j-1}$ is $\omega '$-Giesker 
 stable and $\widetilde{\chi}(F_j \otimes \omega ^{'\otimes n})\ge 
 \widetilde{\chi}(F_{j-1} \otimes \omega ^{'\otimes n})$ for $n\gg 0$.
 Here 
 $\widetilde{\chi}(F\otimes \omega ^{'\otimes n})$ is a 
 reduced Hilbert polynomial of $F\in \Coh (X/S)$. 
 In particular we have $F_j \cong i_{\ast}F_{0,j}$ for some  
 $F_{0,j}\in \Coh (X_0)$ and 
 $\mu _{\omega}(F_{0,j})\ge \mu _{\omega}(F_{0,j-1})$. 
 Note that $F_{0,j}$ is $\omega $-Giesker stable, hence 
 $\mu _{\omega}$-semistable. 
By truncating at $j$ with 
$\mu _{\omega}(F_{j,0})\ge \beta \cdot \omega $, 
 we obtain the exact sequence 
 $$0 \lr E_j \lr E \lr E/E_j \lr 0,$$
 such that $E_j \in \tT '$ and $E/E_j \in \fF '$. $\quad \square$

\vspace{5mm}

Let $\aA _{(\beta ', \omega ')}\subset D(X/S)$
 be the tilting for the torsion pair 
$(\tT ', \fF ')$. 
On the other hand, for
 $\beta '+i\omega ' \in \aA (X/S)_{\mathbb{C}}$, 
define $Z_{(\beta ', \omega ')} \colon \nN(X/S)\to \mathbb{C}$ to 
be 
$$Z_{(\beta ', \omega ')}(E)
 =-\int e^{-(\beta ' +i\omega ')}\ch (E)\sqrt{\td _X}.$$
Under the isomorphism 
$$\Hom (\nN (X/S), \mathbb{C}) \cong \mathbb{C} \oplus 
N^1 (X/S)_{\mathbb{C}} \oplus \mathbb{C},$$
$Z_{(\beta ', \omega ')}$ corresponds to $(1, \beta '+i\omega ', \frac{1}{2}(\beta 
+i\omega )^2)$.  
\begin{lem}
For $\beta ' +i\omega ' \in \aA (X/S)_{\mathbb{C}}$, the function $Z_{(\beta ', \omega ')}$ is a slope function on $\aA _{(\beta ', \omega ')}$ 
if and only if 
for all spherical sheaves $E$ on $X_0$ one has 
$Z_{(\beta ', \omega ')}(i_{\ast}E)\notin \mathbb{R}_{\le 0}$. 
\end{lem}
\textit{Proof}. 
The same proof of~\cite[Proposition 9.2]{Brs2} works. $\quad \square$

\vspace{5mm}

Let us take 
$\beta '+i\omega '\in \aA (X/S)_{\mathbb{Q}}$ 
such that $Z_{(\beta ', \omega ')}(i_{\ast}E)\notin \mathbb{R}_{\le 0}$
for all spherical sheaf $E$ on $X_0$. 
Then since the image of $Z_{(\beta ', \omega ')}$ is 
discrete, the Harder-Narasimhan property is automatically satisfied. 
Thus we obtain the stability condition
$\sigma _{(\beta ', \omega ')}=
(Z_{(\beta ', \omega ')}, \aA _{(\beta ', \omega ')})$ on $D(X/S)$. 
Let $\Stab ^{\circ}(X/S)$ be the connected component of $\Stab (X/S)$ 
which contains $\sigma _{(\beta ', \omega ')}$. 
Applying Theorem~\ref{loc}, we get the following theorem:
\begin{thm}
Let $f\colon X\to S=\Spec \mathbb{C}[[t]]$ be a 
smooth $K3$ fibration. Then 
we have the commutative diagram, 
$$\xymatrix{
\Stab ^{\circ}(X/S) \ar[r]^{\theta ^{\circ}}
 \ar[d]_{\zZ} & \Stab ^{\circ}(X_0) \ar[d]^{\zZ _0} \\
\nN (X)_{\mathbb{C}} \ar[r]^{i^{\ast}} & \nN (X_0)_{\mathbb{C}},
}$$
which gives a homeomorphism between $\Stab ^{\circ}(X/S)$ and 
one of the connected component of $\Stab ^{\circ}(X_0) 
\times _{\nN (X_0)_{\mathbb{C}}} \nN (X)_{\mathbb{C}}$. 
In particular we have the map, 
$$\Stab ^{\circ}(X/S) \lr \pP ^{+}_0(X/S) \cneq \nN (X)_{\mathbb{C}} 
\cap \pP ^{+}_0 (X_0), $$
which is a regular covering map. 
\end{thm}
\textit{Proof}. It is clear that the map $\theta$ constructed 
in the previous section takes $\Stab ^{\circ}(X/S)$ to $\Stab ^{\circ}(X_0)$. 
For $\sigma \in \Stab ^{\circ}(X/S)$ and $U\in \nN (X)_{\mathbb{C}}$, 
 one has 
 $$\lVert U \rVert _{\sigma} = \lVert i^{\ast}U \rVert _{\theta (\sigma)}<
 \infty$$
 by Theorem~\ref{K3}. Thus one can apply 
 Theorem~\ref{loc}. $\quad \square$
 
 \section{(Appendix) Autoequivalences of crepant small resolutions 
 of $cA$-type singularities}
 
 As an appendix, we apply the results in Section 5
 to study the group of autoequivalences for crepant
 small resolutions of $cA$-type singularities. 
 Let $R$ be a three dimensional 
 local complete $\mathbb{C}$-algebra
 which has an isolated $cA_n$-type singularity, i.e. a general hyper plane 
 section $0\in Y\subset \mathrsfs{Y}\cneq \Spec R$ is of $A_n$-type 
 singularity. Let $f\colon X\to Y$ be a minimal resolution 
 and assume that $f$ extends to a crepant small resolution,
 $\widetilde{f}\colon \mathrsfs{X} \to \mathrsfs{Y}$, i.e. one has the 
 Cartesian square:
 $$\xymatrix{
 X \ar[r]^{i} \ar[d]_{f} & \mathrsfs{X} \ar[d]^{\widetilde{f}} \\
 Y \ar[r] & \mathrsfs{Y}, }$$
 and $\widetilde{f}$ is an isomorphism outside $X$. 
 Let $C\subset X$ be 
 the exceptional locus of $f$. Since $Y$ has $A_n$ singularity, 
 $C$ is a chain of rational curves $C=C_1 \cup \cdots \cup C_n$ 
 with $C_i \cap C_j =\emptyset$ for $\lvert i-j \rvert >1$. 
 Let $D(X/Y)$, $D(\mathrsfs{X}/\mathrsfs{Y})$ be the triangulated 
 categories, 
 $$
 D(X/Y)  \cneq \{ E \in D(X) \mid \Supp (E)\subset C \},  \quad
  D(\mathrsfs{X}/\mathrsfs{Y})  \cneq 
 \{ E\in D(\mathrsfs{X}) \mid 
 \Supp (E)\subset C \}, $$
 and denote by $\Stab (\mathrsfs{X}/\mathrsfs{Y})$, $\Stab (X/Y)$ 
 the spaces of locally finite stability conditions on 
 $D(\mathrsfs{X}/\mathrsfs{Y})$, $D(X/Y)$ respectively. 
 Let us define $\Auteq (\mathrsfs{X}/\mathrsfs{Y})$ to be the group 
 $$
 \{ \Phi \colon D(\mathrsfs{X}) \stackrel{\cong}{\lr} D(\mathrsfs{X}) \mid 
 \Phi \mbox{ is of Fourier-Mukai type with kernel supported on }
 \mathrsfs{X}\times _{\mathrsfs{Y}} \mathrsfs{X}\}.$$
 Note that $\Auteq (\mathrsfs{X}/\mathrsfs{Y})$ is 
 regarded as a subgroup of autoequivalences on $D(\mathrsfs{X}/\mathrsfs{Y})$. 
 The similarly defined group $\Auteq (X/Y)$
 was studied in~\cite{UI} and the purpose of this section is to 
 study $\Auteq (\mathrsfs{X}/\mathrsfs{Y})$ via the space 
 $\Stab (\mathrsfs{X}/\mathrsfs{Y})$. 
 Recently the detailed study of the space $\Stab (X/Y)$ has 
 been done 
 by~\cite{IUU}, using the technique of~\cite{UI} and local mirror 
 symmetry. 
 \begin{thm}\emph{\bf{\cite{IUU}}}
 The space $\Stab (X/Y)$ is connected and simply connected.
 \end{thm}
 Then the argument of Section 5 quickly yields the following:
 \begin{thm}\label{connected}
 The space $\Stab (\mathrsfs{X}/\mathrsfs{Y})$ is homeomorphic to 
 $\Stab (X/Y)$. In particular $\Stab (\mathrsfs{X}/\mathrsfs{Y})$
 is connected and simply connected.
 \end{thm}
 On the other hand, we discussed the
  relationship between $\Auteq (\mathrsfs{X}/
 \mathrsfs{Y})$ and $\Stab (\mathrsfs{X}/\mathrsfs{Y})$
 in~\cite{Tst}. We prepare some notations. 
 For $1\le i \le j \le n$, let $C_{i,j}\cneq 
 C_i \cup \cdots \cup C_j \subset C$ and for $k\in \mathbb{Z}$
define $H_{i,j,k}$ to be
$$H_{i,j,k}\cneq \{ \beta +i\omega \in
 N^1 (\mathrsfs{X}/\mathrsfs{Y})_{\mathbb{C}}
\mid (\beta +i\omega)\cdot C_{i,j}=k \}.$$
Let $\Stabn (\mathrsfs{X}/\mathrsfs{Y}) \subset 
\Stab (\mathrsfs{X}/\mathrsfs{Y})$ be the normalized stability 
condition as in Section 4. 
For $\Phi \in \Auteq (\mathrsfs{X}/\mathrsfs{Y})$, we determine
 $n(\Phi)\in \mathbb{Z}$ to be
 $\Phi (\oO _{k(\mathrsfs{X})})=\oO _{k(\mathrsfs{X})}[n(\Phi)]$ for 
 the generic point $\oO _{k(\mathrsfs{X})} \in D(\QCoh(\mathrsfs{X}))$. 
Combing the result of~\cite{Tst} and Theorem~\ref{connected}, one 
obtains the following:
 \begin{thm}\label{crepant}
 One can find a connected component $\Stabn ^{\circ}(\mathrsfs{X}/\mathrsfs{Y})
 \subset \Stabn (\mathrsfs{X}/\mathrsfs{Y})$
 such that $\mathbb{C}\times \Stabn ^{\circ}(\mathrsfs{X}/\mathrsfs{Y})$
 is homeomorphic to $\Stab (\mathrsfs{X}/\mathrsfs{Y})$. 
 In particular $\Stabn ^{\circ}(\mathrsfs{X}/\mathrsfs{Y})$ is 
 simply connected, and one has the map 
 $$\Stabn ^{\circ}(\mathrsfs{X}/\mathrsfs{Y}) \lr 
 N^1 (\mathrsfs{X}/\mathrsfs{Y})_{\mathbb{C}} \setminus \bigcup _{1\le i \le 
 j\le n, k\in \mathbb{Z}}H_{i,j,k},$$
 which provides an universal cover of the right hand side.
 Its Galois group is given by the kernel of the map, 
 $$\Auteq (\mathrsfs{X}/\mathrsfs{Y}) \ni \Phi \longmapsto 
 (\det \Phi (\oO _{\mathrsfs{X}}), n(\Phi)) \in \Pic (\mathrsfs{X})
  \times \mathbb{Z}.$$
 Thus we have 
 $$\Auteq (\mathrsfs{X}/\mathrsfs{Y}) =
 \pi _1 (N^1 (\mathrsfs{X}/\mathrsfs{Y})_{\mathbb{C}} \setminus \bigcup _{1\le i \le 
 j\le n, k\in \mathbb{Z}}H_{i,j,k}) \rtimes \Pic (\mathrsfs{X}) \times
 \mathbb{Z}.$$
 \end{thm}
 Let us investigate $\Auteq (\mathrsfs{X}/\mathrsfs{Y})$ more 
 carefully. Note that for another smooth minimal model 
 $\mathrsfs{W}\to \mathrsfs{Y}$, the induced birational map 
 $\mathrsfs{W}\times _{\mathrsfs{Y}}Y \dashrightarrow X$ 
 extends to an isomorphism. Thus we can regard $C$ as the exceptional 
 locus of $\mathrsfs{W}\to \mathrsfs{Y}$. 
 For $1\le a \le b \le n$, we consider the sequence of 
 birational maps:
 $$
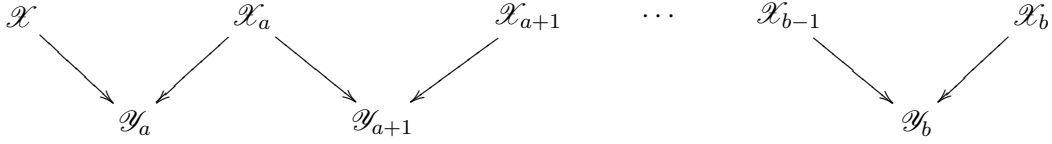
\xymatrix{
 \mathrsfs{X} \ar[dr] &  & \mathrsfs{X}_a \ar[dr] \ar[dl]
 & & \mathrsfs{X}_{a+1} \ar[dl] & \cdots & \mathrsfs{X}_{b-1} 
 \ar[dr]  & & \mathrsfs{X}_b \ar[dl] \\
 & \mathrsfs{Y}_a & & \mathrsfs{Y}_{a+1} & & & & \mathrsfs{Y}_{b} &
 }$$
 Here $\mathrsfs{X}\dashrightarrow \mathrsfs{X}_a$ is a flop at 
 $C_a$ and for each $k$,
 $\mathrsfs{X}_k \dashrightarrow \mathrsfs{X}_{k+1}$ is a 
 flop at $C_{k+1}$. Let $\phi _k \colon \mathrsfs{X}\dashrightarrow 
 \mathrsfs{X}_k$ be the birational map. Let us fix a point 
 $\ast \in \aA (\mathrsfs{X}/\mathrsfs{Y})_{\mathbb{C}}$ and choose 
 a path 
 $$\gamma _{a,b} \colon [0,1] \lr 
 N^1 (\mathrsfs{X}/\mathrsfs{Y})_{\mathbb{C}}\setminus 
 \bigcup _{1\le i\le j\le n, k\in \mathbb{Z}}H_{i,j,k},$$
 which satisfies the following:
 there exists an sequence 
 $$0<t_a <t_{a+1}< \cdots <t_{b-1}<t_b <t_{b+1}=t_{b+1}' <t_b '
 < t_{b-1}' < \cdots < t_a '<1,$$ such that
 \begin{itemize}
 \item $\gamma _{a,b}(0)=\gamma _{a,b}(1)=\ast$. 
 \item $\gamma _{a,b}(0,t_a)$, $\gamma _{a,b}(t_a',1)$ 
 are contained in $\aA (\mathrsfs{X}/\mathrsfs{Y})_{\mathbb{C}}$ and 
 for each $k$, $\gamma _{a,b}(t_k, t_{k+1})$, $\gamma _{a,b}(t_{k+1}', 
 t_k')$ are contained in $\phi _{k\ast}^{-1}\aA (\mathrsfs{X}_k/
 \mathrsfs{Y})_{\mathbb{C}}$. 
 \item For each $k<b-1$, $\gamma _{a,b}(t_{k+1})$, $\gamma_{a,b}(t_{k+1}')$
 are contained in general points of 
 $\phi _{k\ast}^{-1}C_{k+1}(\mathrsfs{X}_k)$ and 
 $\gamma (t_b)$, $\gamma (t_b')$ are contained in general points of 
 $\phi _{b-1 \ast}^{-1}C_b (\mathrsfs{X}_{b-1})$, 
 $\phi _{b-1 \ast}^{-1}C_b'(\mathrsfs{X}_{b-1})$ respectively.
 \end{itemize}
 Here for another model $\mathrsfs{W}\to \mathrsfs{Y}$, we have defined 
 $C_k(\mathrsfs{W})$, $C_k'(\mathrsfs{W})$ to be 
 \begin{align*}
 C_k (\mathrsfs{W}) & \cneq \{ \beta +i\omega \in N^1 (\mathrsfs{W}/
 \mathrsfs{Y})_{\mathbb{C}} \mid -1<\beta \cdot C_k <0, \omega \cdot C_k=0
 \}, \\
 C_k '(\mathrsfs{W}) & \cneq \{ \beta +i\omega \in N^1 (\mathrsfs{W}/
 \mathrsfs{Y})_{\mathbb{C}} \mid 0<\beta \cdot C_k <1, \omega \cdot C_k=0
 \}
 \end{align*}
 In other words, $\gamma _{a,b}$ is an element 
 $$\gamma _{a,b}\in \pi _1 (N^1 (\mathrsfs{X}/\mathrsfs{Y})_{\mathbb{C}}
 \setminus \bigcup _{1\le i\le j\le n, k\in \mathbb{Z}}H_{i,j,k}, \ast),$$
 which goes around $H_{a,b,0}$. Let us describe the autoequivalence of 
 $D(\mathrsfs{X}/\mathrsfs{Y})$ induced by $\gamma _{a,b}$. 
 \begin{defi}
 For $1\le a\le b \le n$, define $E_{a,b}\in \Coh (X)$ as 
 follows: $E_{a,a}=\oO _{C_a}(-1)$ and for $a\le k \le b$, construct
 $E_{a,k}$ successively as the unique non-trivial extension, 
 $$0 \lr E_{a,k-1} \lr E_{a,k}\lr \oO _{C_k}(-1) \lr 0.$$
 \end{defi}
 Note that $E_{a,b}\in \Coh (X)$ is a spherical object
 in $D(X/Y)$. But 
 if we regard $E_{a,b}$ as an object of $D(\mathrsfs{X}/\mathrsfs{Y})$, it 
 is not necessary spherical. 
 Instead, we can use the generalized notion of spherical
 objects and associated twists introduced in~\cite{Tgen}.
 Let us consider the moduli theory of simple sheaves on 
 $\mathrsfs{X}$. Since the dimension of 
 $\Ext _{\mathrsfs{X}}^1 (E_{a,b}, E_{a,b})$ is one or zero, 
 the universal deformation space of $E_{a,b}$ as a sheaf on $\mathrsfs{X}$
 is written as $\Spec \mathbb{C}[t]/(t^{m+1})$ for some $m$. 
 Let us denote $R_m \cneq \mathbb{C}[t]/(t^{m+1})$ and 
 let $\eE _{a,b}\in \Coh (\mathrsfs{X}\times \Spec R_m)$ be the 
 universal family. Then~\cite[Proposition 4.3]{Tgen}
  and~\cite[Remark 4.4]{Tgen} imply $\eE _{a,b}$ is a $R_m$-spherical 
  object in the sense of~\cite[Definition 2.1]{Tgen}.
  Thus one can associate the autoequivalence,
  $T_{\eE_{a,b}} \in \Auteq (\mathrsfs{X}/\mathrsfs{Y})$
  which fits into the triangle~\cite[Theorem 1.1]{Tgen}:
  $$\dR \Hom _{\mathrsfs{X}}(\pi _{\ast}\eE_{a,b}, F)
   \dotimes _{R_m} \pi _{\ast}\eE_{a,b} \lr 
F \lr T_{\eE_{a,b}}(F), $$ 
for $F\in D(\mathrsfs{X})$ and $\pi \colon \mathrsfs{X}\times \Spec R_m \to 
\mathrsfs{X}$ is a projection. 
\begin{lem} \label{loop}
The autoequivalence of $D(\mathrsfs{X}/\mathrsfs{Y})$ induced by 
$\gamma _{a,b}$ coincides with $T_{\eE _{a,b}}$. 
\end{lem}
\textit{Proof}. 
We have a sequence of standard equivalences, 
$$ D(\mathrsfs{X}_b) \stackrel{\Phi _b}{\lr}
D(\mathrsfs{X}_{b-1}) \stackrel{\Phi _{b-1}}{\lr} \cdots 
\stackrel{\Phi _{a-1}}{\lr} D(\mathrsfs{X}_a) 
\stackrel{\Phi _a}{\lr} D(\mathrsfs{X}).$$
Here $\Phi _k$ takes $\iPPer (\mathrsfs{X}_k/\mathrsfs{Y}_k)$ to 
$\oPPer (\mathrsfs{X}_{k-1}/\mathrsfs{Y}_k)$. It is clear from the chamber 
structure on $\Stab (\mathrsfs{X}/\mathrsfs{Y})$ described
 in~\cite[Theorem 4.11]{Tst}
 that the autoequivalence induced by $\gamma _{a,b}$ is given by 
 $$\Phi _{a,b}\cneq 
 \Phi _a \circ \cdots \circ \Phi _{b-1} \circ \Phi _b ^2 \circ 
 \Phi _{b-1}^{-1}\circ \cdots \circ \Phi _a ^{-1}.$$
 On the other hand, let $T_k \cneq T_{\oO _{C_k}(-1)}$ be the usual 
 spherical twist on $D(X)$. Then if we replace $\Phi _k$ by 
 $T_k$, the above composition becomes
 \begin{align*}
 & T_a \circ \cdots \circ T_{b-1}\circ T_b ^2 \circ T_{b-1}^{-1}
 \circ \cdots \circ T_a ^{-1} \\
 &= (T_a \circ \cdots \circ T_{b-1} \circ T_b \circ T_{b-1}^{-1} 
 \circ \cdots T_{a}^{-1})^2 \\
 &= T_{T_a \circ \cdots \circ T_{b-1}(\oO _{C_b}(-1))}^2 \\
 &= T_{E_{a,b}}^2.
 \end{align*}
 Now we have the three commutative diagrams of functors, 
 $$\xymatrix{
 D(\mathrsfs{X}) \ar[r]^{\Phi _{a,b}}  & D(\mathrsfs{X})  \\
 D(X) \ar[u]^{i_{\ast}}\ar[r]_{T_{E_{a,b}}^2} & D(X), \ar[u]_{i_\ast}
 } \qquad
 \xymatrix{
 D(\mathrsfs{X}) \ar[r]^{T_{\eE _{a,b}}}  & D(\mathrsfs{X})  \\
 D(X) \ar[u]^{i_{\ast}}\ar[r]_{T_{E_{a,b}}^2} & D(X), \ar[u]_{i_\ast}
 } \qquad
 \xymatrix{
 D(\mathrsfs{X}) \ar[r]^{T_{\eE _{a,b}}\circ \Phi _{a,b}^{-1}}
   & D(\mathrsfs{X})  \\
 D(X) \ar[u]^{i_{\ast}}\ar[r]_{\id} & D(X). \ar[u]_{i_\ast}
 }$$
 The first diagram follows from Lemma~\ref{fun} below, the 
 second one follows from~\cite[Theorem 4.5]{Tgen}, and the last one 
 follows from the previous two diagrams. 
 Thus for any closed point $x\in \mathrsfs{X}$, 
 $T_{\eE _{a,b}}\circ \Phi _{a,b}^{-1}$ takes $\oO _x$ to $\oO _x$. 
 In this situation, it is well-known that 
 $T_{\eE _{a,b}}\circ \Phi _{a,b}^{-1}=\otimes \lL$
 for some $\lL \in \Pic (\mathrsfs{X})$.
 Again the above commutative diagram implies $\lL =\oO _{\mathrsfs{X}}$, 
 therefore $\Phi _{a,b}=T_{\eE _{a,b}}$. $\quad \square$ 
  
  \hspace{5mm}
  
 We have used the following lemma:
 
 \begin{lem}\label{fun}
 One has the commutative diagram of functors, 
 $$\xymatrix{
 D(\mathrsfs{X}_a) \ar[r]^{\Phi _a}  & D(\mathrsfs{X})  \\
 D(X) \ar[u]^{i_{a\ast}}\ar[r]_{T_a} & D(X). \ar[u]_{i_\ast}
 }$$
 Here we have identified $\mathrsfs{X}_a \times _{\mathrsfs{Y}}Y$ with 
 $X$ and denoted by $i_a$ the inclusions. 
 \end{lem}
 \textit{Proof}. Chen's lemma~\cite[Lemma 6.2]{Ch} yields an equivalence 
 $\Phi _{a,0}\colon D(X) \to D(X)$ which fits into 
 the above commutative diagram after replacing
  $T_{a}$ by $\Phi _{a,0}$. 
  On the other hand
  by~\cite[Proposition 3.5.8]{MVB}, the simple 
  objects of $\iPPer (\mathrsfs{X}_a/\mathrsfs{Y}_a)$,
   $\oPPer (\mathrsfs{X}/\mathrsfs{Y}_a)$ are given by
  $$\{ \oO _{C_a}(-1)[1], \oO _{C_a} \}, \qquad
  \{ \oO _{C_a}(-1), \oO _{C_a}(-2)[1] \}$$ respectively. 
  By~\cite[Lemma 5.1]{Tst}, we have 
  $\Phi _a (\oO _{C_a}(-1)[1])=\oO _{C_a}(-1)$. Therefore 
  $\Phi _a (\oO _{C_a})=\oO _{C_a}(-2)[1]$, and 
  \begin{align*}
  \Phi _{a,0} (\oO _{C_a}(-1)[1])&=\oO _{C_a}(-1)=T_{a}(\oO 
  _{C_a}(-1)[1]), \\
  \Phi _{a,0} (\oO _{C_a})&=\oO _{C_a}(-2)[1]=T_{a}(\oO 
  _{C_a}). \end{align*}
  Thus $\Phi _{a,0}^{-1}\circ T_{a}$ is identity on 
  $\{ \oO _{C_a}(-1), \oO _{C_a} \}$. Since we have the exact sequence 
  $0 \to \oO _{C_a}(-1) \to \oO _{C_a} \to \oO _x \to 0$ for a closed 
  point $x\in C_a$, $\Phi _{a,0}^{-1}\circ T_{a}$
    takes closed points to closed points. Now as in Lemma~\ref{loop}, 
  $\Phi _{a,0}^{-1}\circ T_{a}$ is written as a composition 
  of a pull-back of $\Aut (X)$ and a tensoring a line bundle.
  Therefore $\Phi _{a,0}^{-1}\circ T_{a}$ must be identity
 since it is identity on $\{ \oO _{C_a}(-1), \oO _{C_a} \}$
 and outside $C_a$. $\quad \square$
  
  \hspace{5mm}
  
 Now we can find a generator of $\Auteq (\mathrsfs{X}/\mathrsfs{Y})$. 
 
 \begin{thm} 
 We have 
 $$\Auteq (\mathrsfs{X}/\mathrsfs{Y})=\langle 
 T_{\eE _{a,b}}, \Pic (\mathrsfs{X})
 \mid 1 \le a \le b \le n \rangle \times \mathbb{Z}.$$
 Here $\langle 
 T_{\eE _{a,b}}, \Pic (\mathrsfs{X})
 \mid 1 \le a \le b \le n \rangle$ is a subgroup
 generated by 
 $T_{\eE _{a,b}}$ for $1\le a \le b \le n$ and $\otimes \lL$
 for $\lL \in \Pic (\mathrsfs{X})$, and $\mathbb{Z}$ is generated
 by the shift functor $[1]$. 
 \end{thm}
 \textit{Proof}. 
 The conjugate action of $\Pic (\mathrsfs{X})$ on the 
 subset 
 $$\{ \gamma _{a,b} \mid 1\le a\le b\le n \}
 \subset \pi _1 (N^1 (\mathrsfs{X}/\mathrsfs{Y})_{\mathbb{C}}\setminus 
 \bigcup _{1\le a\le b\le n, k\in \mathbb{Z}}H_{a,b,k},\ast)$$
 provides loops which go around all the codimension two hyper planes
 $H_{a,b,k}$. Thus Theorem~\ref{crepant} and Lemma~\ref{loop}
  imply the result. 
 $\quad \square$

Yukinobu Toda, Graduate School of Mathematical Sciences, University of Tokyo

\textit{E-mail address}:toda@ms.u-tokyo.ac.jp
\end{document}